\documentclass[]{article}
\usepackage[affil-it]{authblk}
\usepackage{graphicx}
\usepackage{fullpage}
\usepackage{amsmath, amssymb,amsfonts}
\usepackage{booktabs}
\usepackage[ruled,vlined,noend]{algorithm2e}
\usepackage{natbib}
 \bibpunct[, ]{(}{)}{,}{a}{}{,}%
\usepackage[dvipsnames]{xcolor}
\usepackage{caption}
\usepackage{subcaption}
\usepackage{comment}
\usepackage{lscape}
\usepackage[normalem]{ulem}

\newcommand{\red}[1]{{\color{black}#1}}

\newcommand{\blue}[1]{{\color{black}#1}}

\newtheorem{definition}{Definition}

\def\GRSC{{GRSC}}
\def\GRSCCB{{GRSC-CB}}
\def\GRSCB{{GRSC-B}}
\def\GRSCC{{GRSC-C}}

\usepackage{hyperref}
\begin{document}
	
	\title{
		The Generalized Reserve Set Covering Problem with Connectivity and Buffer Requirements\footnote{\textcopyright 2019. This manuscript version is made available under the CC-BY-NC-ND 4.0 license \url{http://creativecommons.org/licenses/by-nc-nd/4.0/}.\sloppy Accepted for publication in European Journal of Operational Research; doi: 10.1016/j.ejor.2019.07.017}}
	
	\author[1]{Eduardo \'Alvarez-Miranda \thanks{ealvarez@utalca.cl}}
	
	\author[2]{Marcos Goycoolea
		\thanks{marcos.goycoolea@uai.cl}}
	
	\author[3]{Ivana Ljubi\'c
		\thanks{ivana.ljubic@essec.edu}}
	
	\author[4]{Markus Sinnl\thanks{markus.sinnl@jku.at}}
	
	\affil[1]{Department of Industrial Engineering, Universidad de Talca, Curic\'o, Chile}
	
	\affil[2]{School of Business, Universidad Adolfo Iba\~nez, Santiago, Chile}
	\affil[3]{ESSEC Business School of Paris, France}
	
	\affil[4]{Department of Statistics and Operations Research, Faculty of Business, Economics and Statistics, University of Vienna, Vienna, Austria \newline Institute of Production and Logistics Management, Johannes Kepler University Linz, Linz, Austria}

	\date{}
	\maketitle

\begin{abstract}
The design of nature reserves is becoming, more and more, a crucial task for ensuring the conservation of endangered wildlife. 
In order to guarantee the preservation of species and a general ecological functioning, 
the designed reserves must typically verify a series of spatial requirements. Among the required characteristics, practitioners and researchers have pointed out two crucial aspects: (i) connectivity, so as to avoid spatial fragmentation, and (ii) the design of buffer zones surrounding (or protecting) so-called core areas. 

In this paper, we introduce the \emph{Generalized Reserve Set Covering Problem with Connectivity and Buffer Requirements}. This problem extends the classical \emph{Reserve Set Covering Problem} and allows
to address these two requirements simultaneously.
A solution framework based on Integer Linear Programming and branch-and-cut is developed. The framework is enhanced by valid inequalities, a construction \red{and a primal heuristic and local branching.} \blue{The problem and the framework are presented in a modular way to allow practitioners to select the constraints fitting to their needs and to analyze the effect of e.g., only enforcing connectivity or buffer zones.}

An extensive computational study on grid-graph instances and real-life 
instances based on data from three states of the U.S.
and one region of Australia is carried out to assess the suitability of the proposed model to deal with the challenges faced by decision-makers in natural reserve design. \blue{In the study, we also analyze the effects on the structure of solutions when only enforcing connectivity or buffer zones or just solving a generalized version of the classical reserve set covering problem.} The results show, on the one hand, 
the flexibility of the proposed models to provide solutions according to the decision-makers' requirements,
and on the other hand, the effectiveness of the devised algorithm for providing good solutions in reasonable computing times. 

\end{abstract}

\section{Introduction and Motivation}
\label{sec:intro}

Demographic expansion, natural resource exploitation, and
the consequences of climate change, are among the processes 
that had resulted in a dramatic loss of biodiversity in the last decades. 
In July 2012, at the Rio+20 Earth Summit,
the International Union for Conservation of Nature revealed that about 20,000 species  are threatened with extinction. 
Among them, 4,000 are described as \emph{critically endangered} and 6,000 as \emph{endangered}, while more than 10,000 species are listed as \emph{vulnerable}~\citep[][]{IUCN}.

The maintenance of biodiversity is crucial for the humankind, 
so its preservation is decisive for future generations~\citep[][]{CardinaleEtAl2012}. 
As a matter of fact, immense efforts  have been devoted in the last decades by international organizations, governments, and foundations, for the establishment of protected areas aiming at ensuring a sustainable landscape for wildlife.
The reader is referred to the books~\citep[][]{GergelTurner2002,LindenmayerFranklin2002,MillspaughThompson2008} for in-depth analyses and discussion of motivations, models, cases, and challenges in the field of wildlife conservation and nature reserve planning.

The design of such nature reserves has lead to the development of a plethora of mathematical models. The purpose of such models is in providing optimally designed reserves that respect 
ecological, economical and eventually other requirements (see \cite{{Billionnet2013}} for a recent comprehensive review on optimization models for biodiversity conservation).  
In its most fundamental form, a nature reserve design problem can be stated as follows: 
one is given a set $V$ of land sites  (also known as land units or parcels)  eligible for selection, a set of species or features $S$, and, 
for each species $s \in S$, a set of suitable land sites $V_s\subseteq V$. 
The goal is to find  \red{the} \emph{minimum} number of
reserve sites such that each species is present in the selected set of sites at least once.
This problem is known as the \emph{Reserve Set Covering Problem} (RSC)~\citep[][]{ChurchEtAl1996,PresseyEtAl1997}. 

The RSC can be formulated as an Integer Linear Programming (ILP) problem; 
let $\mathbf{x}\in\{0,1\}^{|V|}$ be a vector of binary variables such that $x_i = 1$ if site $i\in V$ is selected, and $x_i = 0$ otherwise. 
Using this notation, the model 
\begin{align}
\text{(RSC)}\quad\quad\min \pi(\mathbf{x})= \sum_{i\in V} x_ i&\label{eq:RSC1}\tag{CARD}\\
\text{s.t.}\quad \sum_{i\in V_s} x_i &\geq 1,\;\forall s\in S\label{eq:RSC2}\tag{COV}\\
\mathbf{x}&\in\{0,1\}^{|V|},\label{eq:RSC3} \tag{BIN-$\mathbf{x}$}
\end{align}
allows to find the optimal (minimum size) reserve given by $V^* = \{i\in V\mid x_i^*=1\}$. 
Observe that in the RSC, we are only given a set of land sites $V$. 
To work with spatial requirements, such as connectivity and buffer zones, 
one is also given a graph $G=(V,E)$, 
in which the set of edges $E$ encodes the spatial relationship of the land sites (e.g., $\{i,j\} \in E$, if and only if $i$ and $j$ from $V$ share a common border). 
\red{In} the remainder of the paper, the terms \emph{land site} and \emph{node} will be used interchangeably.

The use of mathematical optimization models, 
such as~\eqref{eq:RSC1}-\eqref{eq:RSC3}, 
for the optimal design of nature reserves is a widely explored research area; 
some references covering two decades of work are found in~\citep[][]{Beyer201614,ConservationCorrTech, ClemensEtAl1999,OhmanLamas2005,OnalBriers2003,Williams2008,WilliamsRevelle1996,WilliamsReVelle1998,WilliamsEtAl2004}.
A natural alternative to this problem, is to find a set of exactly $p$ parcels that \emph{maximize} the number of protected species; this problem is known as the Maximal Covering Species Problem~\citep[][]{ChurchEtAl1996}.

Although the RSC enables decision makers to gain important insights about the suitable sites that need to be preserved, 
the obtained solutions typically fail in satisfying relevant 
spatial attributes. 
According to the reviews presented in~\citep[][]{Billionnet2013} and~\citep[][]{WilliamsEtAl2005}, 
the spatial requirements that are typically imposed when designing reserves can be classified as follows: (i) reserve size or \emph{compactness}, (ii) \emph{number} of reserves, (iii) reserves \emph{proximity}, (iv) reserve \emph{connectivity}, (v) reserve \emph{shape}, and (iv) presence of \emph{core} and \emph{buffer} areas.
As pointed out in the above mentioned reviews (and the references therein), 
these requirements correspond to different ecological needs which depend on the considered landscape, 
the species to be protected, and the human activities surrounding the potential reserve. 
Fundamental works on the design of reserves respecting spatial requirements can be found, for example, in~\citep[][]{OnalBriers2002,Schwartz1999}.

Among relevant spacial requirements, connectivity is one of the most prevalent ones; it avoids habitat fragmentation improving the conditions for sustainable ecosystems~\citep[][]{BeierNoss1998,DebinskiHolt2000,GergelTurner2002,MillspaughThompson2008}. 
Various modeling approaches have been proposed for incorporating connectivity into the design of nature reserves \blue{(see,~e.g.,~\citep{Beyer201614,Billionnet2012,billionnet2016designing,briers2002incorporating,cerdeira2010species, DilkinaGomes2010,jafari2013new,jafari2017achieving,OhmanLamas2005,onal2006optimal,onal2016optimal,StJohnTothZabinsky2018,urban2001landscape,WangOnal2011,wang2013designing,williams2002zero}).} 
Note that ensuring connectivity \red{does not necessarily mean that only a single reserve has to be designed. Instead, multiple connected components may be allowed too.}
%
%
\red{The second important criterion is } the presence of core areas and buffer zones that allows the development of so-called \emph{Biosphere} reserves~\citep[][]{Batisse1982,Batisse1990}.
The purpose of the buffer zone is to surround the
core area and protect it from negative external impacts, therefore promoting the long-term viability of critical species (see,~e.g., Chapters~1 and~5 of \cite{LindenmayerFranklin2002} and Chapters~19 and~20 of \cite{MillspaughThompson2008}).
Finally, another important issue in reserve design is to ensure \red{\emph{minimum 
quotas of ecological suitability}} for some species, 
especially critically endangered ones (see, e.g., Chapters~9 and~14 of~\cite{MillspaughThompson2008}).

\paragraph{Contribution and Outline of the Paper}
Although the design of nature reserves considering buffer zone requirements has been addressed before 
(see, e.g.,~\cite{ClemensEtAl1999,WilliamsRevelle1996,WilliamsReVelle1998}), 
the question on how to impose connectivity requirements to the buffer zones remained unanswered in the existing literature.
In particular, in~\citep[][]{Billionnet2013}, the authors point out that known models fail in providing suitable conditions for particular endangered species precisely due to the lack of connectivity of the resulting reserves. 
The main contribution of our work consists \red{thus} of providing, for the first time, a modeling and algorithmic framework for the optimal design of wildlife reserves \red{by simultaneously integrating three criteria}: 
connectivity requirements, construction of buffer zones \red{and minimum 
quotas of ecological suitability}.
This is done by introducing the \emph{Generalized Reserve Set Covering Problem with Connectivity and Buffer Requirements (\GRSCCB)}.
The \GRSCCB\ allows to design a reserve
comprised by one or more connected components; each of them consisting of a
core surrounded by a buffer zone. \blue{The \GRSCCB\ is introduced in a modular way, i.e., we first introduce a generalization of the RSC denoted as \emph{Generalized Reserve Set Covering Problem}, followed by the \emph{Generalized Reserve Set Covering Problem with Buffer Requirements} and \emph{Generalized Reserve Set Covering Problem with Connectivity Requirements}, before putting all constraints together to obtain the \GRSCCB}.
In our extensive numerical study on grid-graph and real-life instances, 
we demonstrate that the \GRSCCB, and its variants,
deliver solutions that properly embody different spatial and ecological features and significantly improve upon the spacial structure of reserves created by the simpler RSC models from the literature. \blue{Moreover, the intermediate models introduced in our work allow decision makers to analyse tradeoffs involved with requiring connectivity or buffer zones.}

The article is organized as follows.
In Section~\ref{sec:MIPFormulation} we incrementally develop an ILP-formulation for the \GRSCCB\ by first giving a generalization of the RSC without considering connectivity and buffer requirements, and then 
introducing constraints to model these two requirements. 
In Section~\ref{sec:algDesc}, we describe a branch-and-cut framework to solve the proposed formulation. 
Computational results on synthetic instances, as well as case studies on real-life instances based on data of the U.S. National Gap Analysis Program
and of the Northern Australia Water Futures Assessment Program are presented in Section~\ref{sec:compres}. Section~\ref{sec:conclu} concludes the paper.

\section{The Generalized Reserve Set Covering Problem with Connectivity and Buffer Requirements}
\label{sec:MIPFormulation}

\begin{figure}[h!tb]
\centering
\begin{subfigure}[b]{0.23\textwidth}
  \includegraphics[scale =0.3] {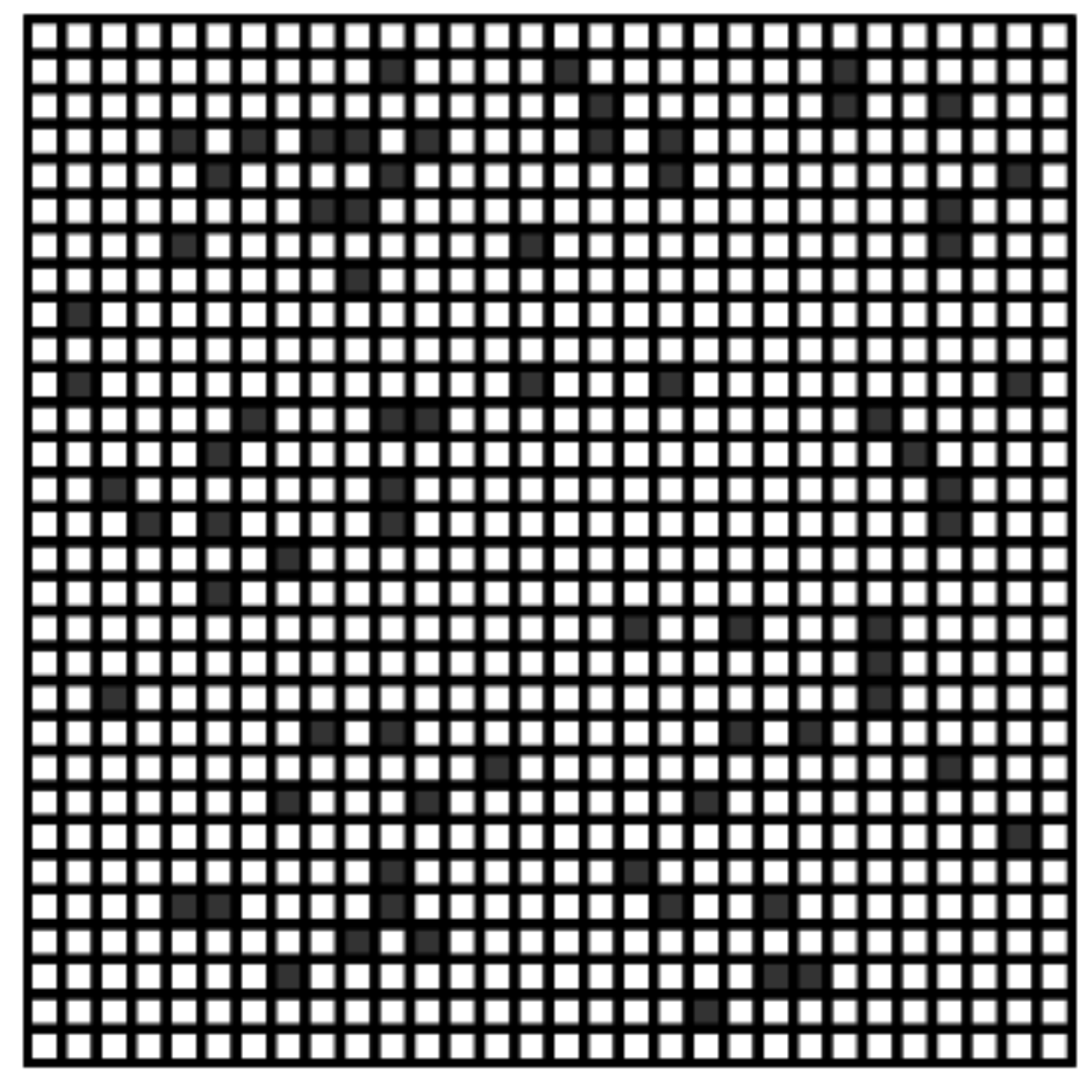}
  \caption{\GRSC\ Solution}
   \label{fig:RSC}
  \end{subfigure}
  \hfill
\begin{subfigure}[b]{0.23\textwidth}
  \includegraphics[scale =0.3] {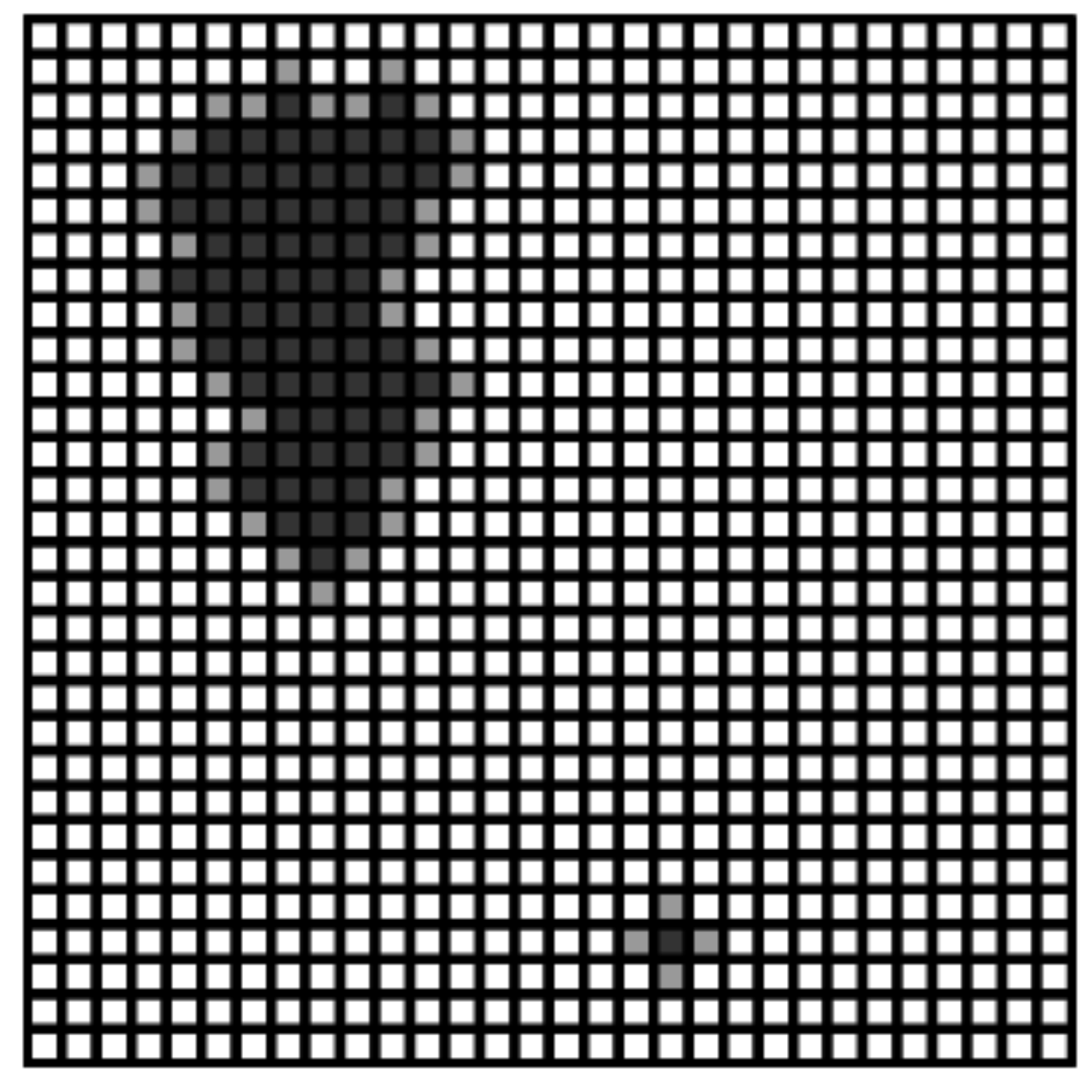}
  \caption{\GRSCB\ Solution}
   \label{fig:ConRSC}
 \end{subfigure}
   \hfill
\begin{subfigure}[b]{0.23\textwidth}
   \includegraphics[scale =0.3] {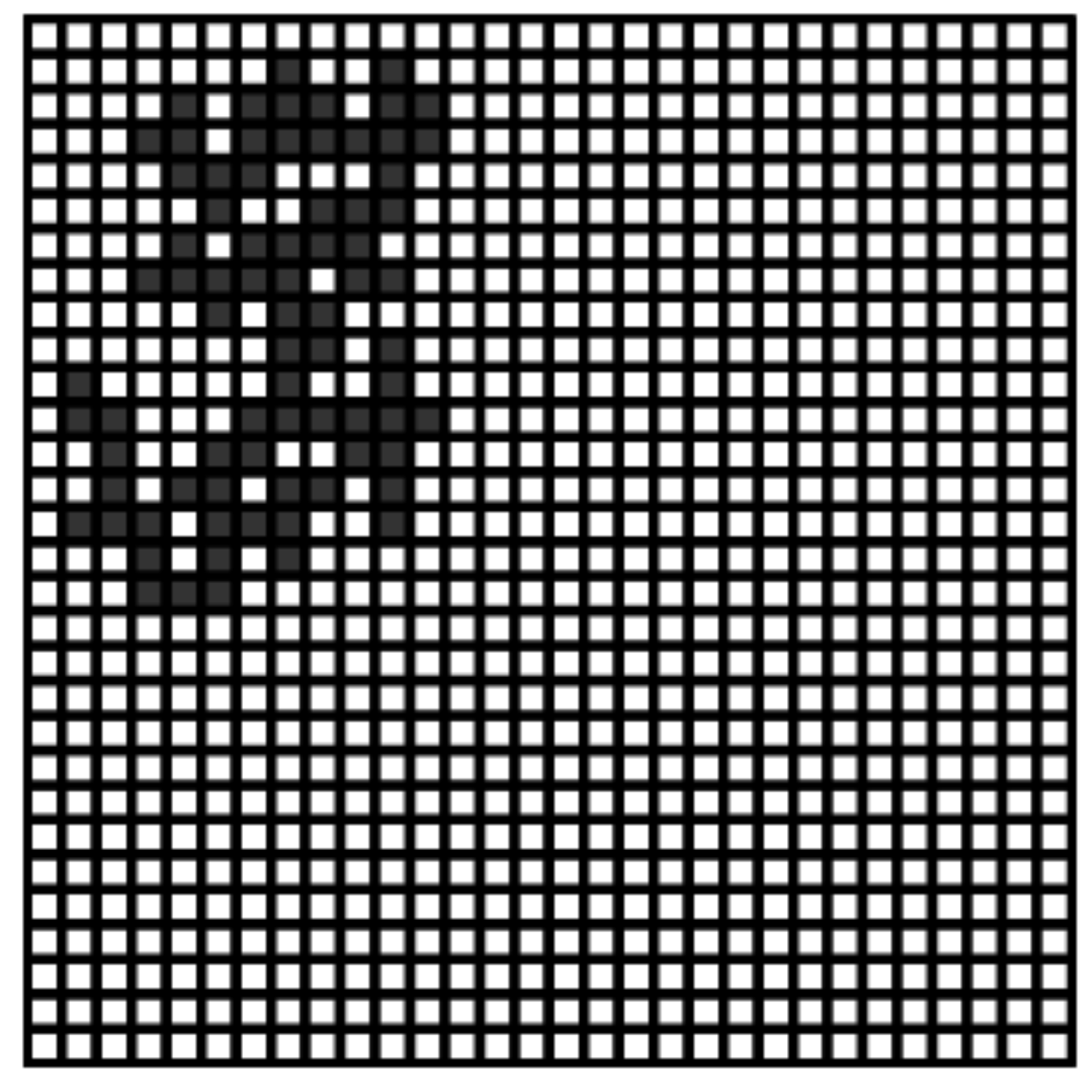}
   \caption{\GRSCC\ Solution}
    \label{fig:BuffRSC}
 \end{subfigure}
   \hfill
\begin{subfigure}[b]{0.23\textwidth}
     \includegraphics[scale =0.3] {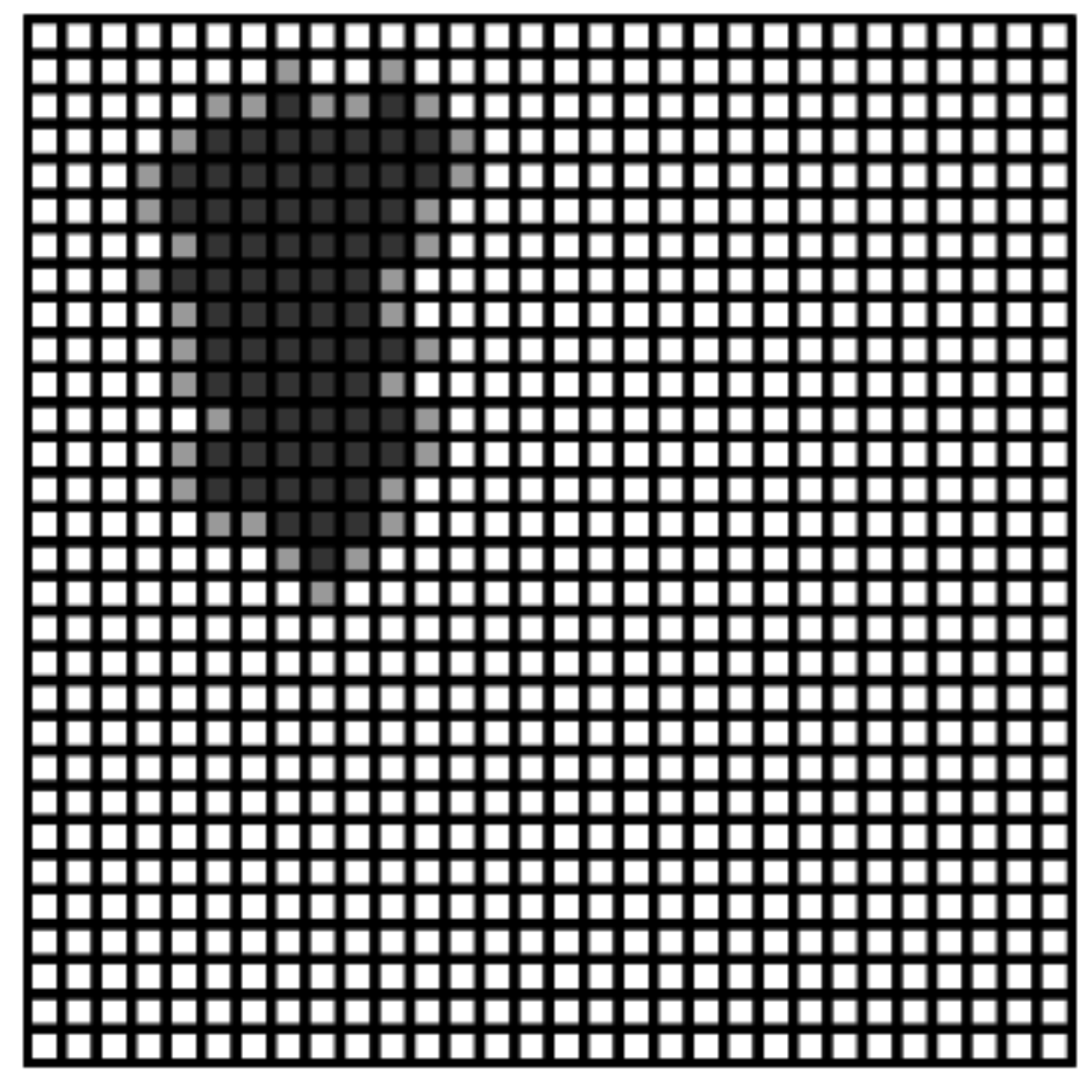}
     \caption{\GRSCCB\ Solution}
      \label{fig:ConBuffRSC}
 \end{subfigure}    
\caption[]{{{
Optimal solutions on the same instance for the different problem variants introduced in the paper.
For problems with buffer, i.e., \GRSCB\ and \GRSCCB,
\protect \raisebox{0.3em}{\fcolorbox{black}{black!70!white}{\null}} denotes core land sites, and \raisebox{0.3em}{\fcolorbox{black}{black!35!white}{\null}}} denotes buffer land sites}; for the problem without buffer, i.e., \GRSC\ and \GRSCC, \protect \raisebox{0.3em}{\fcolorbox{black}{black!70!white}{\null}} denotes selected land sites.}
 \label{fig:RSCExamples}
\end{figure}

In this section, we first give a generalization of the RSC, which takes into account that in real-life situations, 
different types of species need different levels of protection. 
The resulting problem is denoted as \emph{Generalized Reserve Set Covering Problem (\GRSC)}. 
We then introduce constraints for modeling buffer and connectivity requirements, and finally obtain the \GRSCCB. 
These constraints can also be added individually to obtain problems that we denote as \emph{\GRSC\ with buffer requirements (\GRSCB)} and \emph{\GRSC\ with connectivity requirements (\GRSCC)}.
In Figure~\ref{fig:RSCExamples} optimal solutions for these four problems on the same underlying instance are shown. The influence of the spatial constraints imposed in the different problem variants can be easily seen. 
The \GRSC\ solution is very fragmented. 
The \GRSCB\ solution has a more compact shape, but consists of two components, and one of these components is very small.
The \GRSCC\ solution is connected, but its spatial distribution lacks 
 compactness.  Finally, the \GRSCCB\ solution has a similar compact shape as the \GRSCB\ solution, but consists of a single component only.

To account for the difference 
\red{concerning the} need of protection among species, we
consider a partition of the set of species $S$ into two subsets, 
namely $S_1$ and $S_2$, $S = S_1 \cup S_2, S_1 \cap S_2 = \emptyset$. 
Species in $S_1$ are the ones needing stronger protection 
(e.g., the land sites selected for hosting them need to be in the \emph{core} of the designed reserve, if a model with core and buffer zones is considered). 

Additionally, ecologists usually do not only know in which land sites a species lives, but they are able to assess the habitat features of the different land sites.
We model this by using a habitat \emph{suitability} score function $\mathbf{w}:V\times S\rightarrow\mathbb{R}_{\geq 0}$, such that $w_i^s$ measures how advisable, with respect to species $s\in S$, it is to select the land site $i\in V$ as part of the reserve~(see,~e.g., Chapters~9 and~10 of~\cite{MillspaughThompson2008}).
A species $s\in S$ is considered to be protected by the designed reserve, 
if the $w_i^s$ values of the land parcels contained in the reserve sum up to at least $\lambda_s \geq 0$, which represents a \emph{minimum quota of ecological suitability} for species $s$. 
If a model with buffer requirements is used, for species $s \in S_1$, 
the $w_i^s$ values of land parcels selected for the core of the designed reserve must sum up to at least $\lambda_s$ 
(as only core area is suitable to protect the species in $S_1$). 
We denote with $V_s$ the set of land parcels with $w^s_i>0$.
We observe that a simpler version of this suitability approach has already been used previously, e.g., in~\cite{Polasky01022001} the constraint~\eqref{eq:RSC2} is replaced by $\sum_{i\in V_s} x_i \geq h_s,\;\forall s\in S$, 
where $h_s\in\mathbb{Z}_{\geq 1}$ is the number of sites required by species $s\in S$.
Finally, there is a cost function $\mathbf{c}:V\rightarrow\mathbb{R}_{> 0}$, 
such that $c_i$ corresponds to the cost of selecting the land site $i\in V$ as part of the reserve (see,~e.g., Chapter~5 of~\cite{MillspaughThompson2008} for further details regarding economic considerations for planning wildlife conservation).

In the typical reserve design setting considered in the literature, 
all species under consideration must be hosted within the reserve.
However, if the area under consideration hosts many different species, the obtained reserves might be impractical due to size considerations or undesirable shapes. 
To allow for more flexibility in the design, we introduce two parameters, $0 \leq P_1\leq |S_1|$ and $0 \leq P_2\leq |S_2|$, allowing the decision maker to specify for how many of the species in $S_1$ and $S_2$ the obtained reserve must fulfill the minimum quota of ecological suitability.
We observe that the classical RSC described in the introduction can be obtained by
setting $S_2=S$, $w=\mathbf 1$ and $\lambda_s=1$ for every species $s$, and $P_2=|S_2|$.

\begin{definition}
Given the input data described above, 
the \GRSCCB\ is the problem of finding a minimum cost natural reserve
that fulfills the following criteria: 
\begin{description}
\item[(i)] at least $P_1$ species from $S_1$ are hosted in the core area, 
\item[(ii)] at least $P_2$ species from $S_2$ are hosted in the reserve, 
\item[(iii)] the minimum ecological suitability quotas for each of the hosted species are satisfied, 
\item[(iv)] the solution is comprised by at most $k$ connected areas, and
\item[(v)] each connected area consists of a core surrounded by a buffer of width $d$.
\end{description}
\end{definition}

In this article we also address three relaxations of this problem, namely: \GRSC, in which the conditions (iv) and (v) are relaxed; \GRSCB, in which condition (iv) is relaxed, and \GRSCC, in which condition (v) is relaxed (i.e., $d=0$). 

In the following, we introduce the ILP models for all four variants in a modular fashion. \blue{The models for problems without connectivity requirements (i.e., \GRSC\ and \GRSCB) are \emph{compact} models, while the models for problems with connectivity requirements (i.e., \GRSCC\ and \GRSCCB) have an exponential number of constraints. These constraints are separated on-the-fly using branch-and-cut, more details on the separation are given in Section \ref{sec:algDesc}.}

\subsection{Modeling the \GRSC\ Problem} 
In our formulation, vector $\mathbf{x}\in\{0,1\}^{|V|}$ is associated with the decision of selecting sites as part of the reserve (either as part of the core or part of the buffer zone). 
In addition, let vector $\mathbf{z}\in\{0,1\}^{|V|}$ be associated with the decision of selecting sites as part of the core area. 
Let $\mathbf{u}\in\{0,1\}^{|S|}$ be a vector of variables
so that $u_s = 1$ if species $s\in S$ is hosted by the reserve, and
$u_s = 0$, otherwise.
A triplet $(\mathbf{u},\mathbf{x},\mathbf{z})$  is associated with an appropriate territorial coverage of the species if the following inequalities hold:
\begin{align}
\sum_{i\in V_s} w_i^s z_i &\geq \lambda_s u_s ,\;\forall s\in S_1\label{eq:coverZ}\tag{$S_1$-SQ}\\
\sum_{i\in V_s} w_i^s x_i &\geq \lambda_s u_s,\;\forall s\in S_2\label{eq:coverX}\tag{$S_2$-SQ}\\
\sum_{s\in S_1} u_s &\geq P_1\label{eq:sizeC}\tag{$S_1$-PROTECT}\\
\sum_{s\in S_2} u_s &\geq P_2.\label{eq:sizeV}\tag{$S_2$-PROTECT}
\end{align}Constraints~\eqref{eq:coverZ} ensure that if a species $s\in S_1$ is hosted by the reserve ($u_s = 1$), 
then the ecological suitability of the core of the reserve with respect to $s$ must be at least $\lambda_s$.
Likewise, constraints~\eqref{eq:coverX} ensure that if a species $s\in S_2$ is
protected ($u_s = 1$), 
then the ecological suitability of the reserve must be at least $\lambda_s$.
These two set of constraints will be referred to as \emph{suitability quota constraints}. 
Note that if a planner wants to ensure that a given species $s$ \emph{must}
be part of the reserve, then she/he can achieve this by simply adding the constraint $u_s = 1$. 
Constraints~\eqref{eq:sizeC} and~\eqref{eq:sizeV} imposes that at least $P_1$ core species and $P_2$ buffer species must be protected. 
These two constraints will be referred to as \emph{species protection constraints}.

The cost of the reserve is calculated as the sum of the cost of all of the selected sites. 
To model this, we need the following set of linking constraints
\begin{align}
z_i \leq x_i,\;\forall i\in V, \label{eq:linking}\tag{LINK}
\end{align}
i.e., if a site is considered to be in the core, the site must be in the designed reserve. 
The cost of a reserve encoded by $(\mathbf{u},\mathbf{x},\mathbf{z})$ is then given by
\begin{align}
\gamma(\mathbf{u},\mathbf{x},\mathbf{z}) = \sum_{i\in V} c_i x_i.\label{eq:cost}\tag{COST}
\end{align}

The inequalities presented so far allow to formulate a first generalization of the RSC, denoted as \emph{Generalized Reserve Set Covering Problem (\GRSC)}, given by
{\small\begin{align}
(\mbox{\GRSC})\quad\min \bigg\{\gamma(\mathbf{u},\mathbf{x},\mathbf{z})\bigg| \eqref{eq:coverZ},\eqref{eq:coverX}, \eqref{eq:sizeC},\eqref{eq:sizeV}, \eqref{eq:linking},(\mathbf{u},\mathbf{x},\mathbf{z})\in\{0,1\}^{|S|+2|V|}\bigg\}.\nonumber
\end{align}}
The \GRSC\ \red{takes into account the minimum suitability quotas, but it } ensures neither the connectivity of the reserve,
nor the existence of a buffer around the core.
This means that in particular situations, like the one depicted in 
Figure~\ref{fig:RSC}, very fragmented reserves could be obtained.

In the following, these two missing properties are characterized and, 
thereafter, the complete ILP formulation for \GRSCCB\ is presented. 
(We observe that the \GRSC\ could be modeled only using $(\mathbf{x,u})$;
$\mathbf{z}$ will be used in the following for modeling the core/buffer interplay.)

\subsection{Modeling the Buffer Zone} 
To model the buffer surrounding the selected core areas, the following definition of \emph{$d$-neighborhood set of a node $i$} is needed.
\begin{definition} ($d$-neighborhood) For a given integer $d\geq 0$ and a given land site $i\in V$, the $d$-neighborhood of $i$, $\delta_d(i)$,
is defined as
$$
\delta_d(i) = \left\{j\in V_{i\neq j}\mid \mbox{the minimum number of hops between $i$ and $j$ is at most $d$}\right\}.
$$
\red{We also define the set $\delta_d^+(i) = \delta_d(i) \cup \{i\}$.}\end{definition}
Set $\delta_d(i)$ corresponds to the set of all nodes $j\in V$ that are separated by at most $d$ edges from $i$, e.g., all adjacent nodes of $i$ are given by $\delta_1(i)$ (for sake of readability, we will write $\delta_1(i)=\delta(i)$ in the following). As the graph in our setting is undirected, for two nodes $i,j$, we have $i \in \delta_d(j)$ if and only if $j \in \delta_d(i)$.
The definition allows to model the buffer zones with the following set of constraints
\begin{align}
z_i &\leq x_j,\;\forall j\in\delta_{d}(i),\;\forall i\in V\label{eq:buff}\tag{$d$-BUFF.1}
\end{align}
i.e., if $i$ is taken as part of the core ($z_i = 1$), 
then all of the land sites $j\in \delta_d(i)$ must verify $x_j = 1$ 
(i.e., they must at least be part of the buffer zone). 
Hence, it ensures that the each core land site is surrounded by \emph{at least} $d$ other sites within the reserve.
However, using only constraints~\eqref{eq:buff}, it is possible, that some land parcel $i$ is selected (i.e., $x_i=1$ and $z_i=0$) without any parcel $j$ forcing it to be in the solution as a buffer side for $j$. 
This can happen in presence of constraints~\eqref{eq:coverX} and~\eqref{eq:sizeV} (without these constraints, the objective function together with $c_i>0$, $i \in V$ takes care of the issue). 
This situation is not desired, as it can lead to a fragmented reserve. 
To deal with this issue, we propose the following set of inequalities;
\begin{align}
x_i &\leq \sum_{j\in\delta_{d}(i)} z_j,\;\forall i\in V\label{eq:buff2}\tag{$d$-BUFF.2}.
\end{align}
Inequalities~\eqref{eq:buff2} ensure that whenever some $x_i=1$, at least one core node $z_j$ is selected in the $d$-neighborhood. Therefore, combining~\eqref{eq:buff} and~\eqref{eq:buff2} leads to cores that are \emph{nested} within buffer zones of width $d$.
Adding constraints~\eqref{eq:buff} and~\eqref{eq:buff2} to \GRSC\ leads to the \emph{\GRSC\ with Buffer Requirements (\GRSCB)}.

\subsection{Modeling Connectivity of the Reserve} 
\label{subsec:connectivity}


Connectivity is modeled using a variant of a concept called \emph{node-separators} (see,~e.g.,~\cite{AlvarezLjubicMutzel2013a,AlvarezLjubicMutzel2013b,CarvajalEtAl2011,FischettiEtAl2014}).
For modeling purposes, let $G_r = (V_r,E_r)$ with
$V_r = V\cup\{r\}$ and $E_r = E\cup\{(r,i)\mid i \in V \}$, i.e., $G_r$ corresponds to the original graph 
\red{extended by} 
an artificial root $r$ and $|V|$ additional directed arcs connecting $r$ with every node in $V$ (denoted $r$-arcs, $A_r$).
For a given site  $\ell$ from $V$,
a tuple $W = (W_V,W_A)$, $W_V\subseteq V$ and $W_A\subseteq A_r$,
is called an \emph{$r$-arc-node-separator} if and only if after removing
$W$ from $G_r$, site $\ell$ cannot be reached from $r$.
For a given node $\ell$ from $V$,
let $\mathcal{W}_\ell$ be the set of all $r$-arc-node-separators
with respect to $\ell$.

Let $\mathbf{y}\in\{0,1\}^{|V|}$ be a vector of auxiliary variables, such that $y_{i} = 1$ if $i\in V$ is connected with $r$ through arc $(r,i)$,
and $y_{i} = 0$ otherwise. 
Using these variables, we model connectivity (of both the core area and the buffer area) by ensuring that in the obtained solution, 
there is a path from $r$ 
to all core nodes and all buffer nodes selected in the solution (i.e., nodes with $z_i = 1$, resp., $x_i=1$).
This means, the obtained feasible solutions are comprised of connected components rooted at nodes with $y_i=1$.
Thus, the number of connected components within the reserve can be
controlled with the constraint
\begin{align}
\sum_{j\in V} y_j & \leq k.\label{eq:root}\tag{NCOMP}
\end{align}
Constraint~\eqref{eq:root} can also be written in equality form, if desired by a decision maker. Moreover, if a land parcel is connected to the root (i.e., $y_i=1$) the land parcel must be taken in the core; this ensures that every connected component in the solution has at least one core land parcel. This linking is done with the following set of constraints:
\begin{align}
y_i & \leq \red{z_i,}\;\forall i \in V.\label{eq:link}\tag{YZ}
\end{align}

The connectivity inequalities follow easily from the definition of the $r$-arc-node-separators: If a node $\ell$ is in the solution, then for each separator in $W=(W_V,W_A) \in \mathcal W_\ell$, at least one element (i.e., node from $W_V$ or arc from $W_A$) must be taken.
We obtain the following family of inequalities 
for connectivity of the selected core sites
\begin{align}
\sum_{i\in W_V} z_i  +  \sum_{j\in W_{A}} y_{j} &\geq z_\ell,\;\forall  W \in {\cal W}_\ell,\;\forall \ell \in V;\label{eq:concore}\tag{CORECON}
\end{align}
inequalities~\eqref{eq:concore} are denoted as \emph{connectivity cuts}. 
They are exponential in number,
and the resulting ILP is tackled by means of branch-and-cut (see Section~\ref{sec:algDesc} for a separation procedure for these inequalities). 
We observe that these cuts can be down-lifted by not considering all $j \in W_A$, but only $j \leq \ell$, i.e., if node $\ell$ is in the solution, the component containing $\ell$ must be rooted either at $\ell$, or a node with index smaller than $\ell$.
By using this down-lifted variant, all $k$ components of a feasible solution must be rooted at the node with minimal index within each component, e.g., symmetric solutions giving the same components, but rooted at other nodes, are not allowed.

When using inequalities~\eqref{eq:buff} and~\eqref{eq:buff2} together with the inequalities proposed in this section, 
connectivity of the land parcels 
selected by $x_i=1$ is automatically ensured, since all parcels with $x_i=1$ and $z_i=0$ form a buffer of thickness $d$ around all $i$ with $z_i=1$. 
The reserve induced by $z_i=1$ form (at most) $k$ connected components and the buffer around them is connected by construction (as argued in the previous section).

If a decision maker does not want a buffer-zone, but only a reserve comprised of (at most) $k$ connected components, connectivity cuts must be written in $x$ instead of $z$, i.e., 
\begin{align}
\sum_{i\in W_V} x_i  +  \sum_{j\in W_{A}} y_{j} &\geq x_\ell,\;\forall  W \in {\cal W}_\ell,\;\forall \ell \in V;\label{eq:conbuffer}\tag{ALLCON}
\end{align}
and also the \red{linking-constraints~\eqref{eq:linking}}
should be replaced by
\begin{align}
y_i & \leq x_i,\;\forall i \in V.\label{eq:link2}\tag{YX}
\end{align}We define the \emph{\GRSC\ with Connectivity Requirements (\GRSCC)} as the problem obtained by adding~\eqref{eq:conbuffer}, \eqref{eq:link2} and~\eqref{eq:root} to \GRSC.

\subsection{The Complete Model} 

We are now ready to give the ILP model for the \emph{Generalized Reserve Set  Covering Problem with Connectivity and Buffer Requirements (\GRSCCB)}: 
\begin{align}
(\mbox{\GRSCCB})\quad\quad\quad\min \bigg\{\gamma(\mathbf{u},\mathbf{x},\mathbf{z})\bigg| \mbox{\eqref{eq:coverZ},\eqref{eq:coverX}, \eqref{eq:sizeC},\eqref{eq:sizeV}, \eqref{eq:linking}, }
  \mspace{50mu}
  \notag\\
  \mbox{\eqref{eq:buff},\eqref{eq:buff2},\eqref{eq:concore},\eqref{eq:link},\eqref{eq:root}},(\mathbf{u},\mathbf{x},\mathbf{z},\mathbf{y})\in\{0,1\}^{|S|+3|V|}\bigg\}.\nonumber
\end{align}
\red{In order to solve (\GRSCCB), we propose in the following an algorithmic framework based on branch-and-cut. Besides dealing with the separation of inequalities which are exponential in size, this framework also integrates multiple heuristics that allow to find high-quality solutions for large-scale instances, for which the optimal solution cannot be obtained.}

\section{A Branch-and-Cut Framework for the \GRSCCB}
\label{sec:algDesc}

\red{The core of our branch-and-cut framework} is based on the separation of the connectivity cuts~\eqref{eq:concore}. \red{In addition, to improve the quality of lower bounds, we propose additional valid inequalities for the model (\GRSCCB), denoted as \emph{species cuts, cover inequalities, and species-cover cuts}. These inequalities (which are described in Section~\ref{subsec:valineq}) are dynamically  separated, together with~\eqref{eq:concore}, and the underlying separation routines are given in Section~\ref{subsec:sep}.}
In order to find feasible solutions within the framework, 
we also designed heuristics that are detailed in Section~\ref{subsec:heur}.

\subsection{Valid Inequalities}
\label{subsec:valineq}

\paragraph{Species Cuts}
If a species $s \in S_1$ is in the solution, i.e., $u_s=1$, then there must be a path from the root node to at least one land parcel in $V_s$ (which are the nodes with $w^s_i>0$).
To enforce this, we again use connectivity cuts based on $r$-arc-node separators. The cuts are defined using a graph $G_r^s$, which is an extension of $G_r$ with an additional node $s$.
Let $G^s_r=(V^s_r,E^s_r)$ with $V^s_r=V_r \cup \{s\}$ and $E^s_r = E_s\cup\{(i,s)\mid i \in V_s \}$, i.e., $G^s_r$ corresponds to $G_r$ plus
an artificial sink node $s$ and $|V_s|$ additional directed arcs connecting every node in $V_s$ with $s$ (denoted as $A_s$).
For a given species $s \in S_1$, let $\mathcal{W}_s$ be the set of all $r$-arc-node-separator with respect to $s$. 
The species cuts for a given species $s$ read as follows
\begin{align}
\sum_{i\in W_V} z_i  +  \sum_{j\in W_A} y_j \geq u_s,\;\forall W \in \mathcal W_s, s\in S_1.\label{eq:zgSC}\tag{SC}
\end{align}

\paragraph{Cover Inequalities and Species-cover Cuts} 
Let $W_s=\sum_{i \in V} w^s_i $. For a given species $s\in S$, a set $C_s\subset V_s$ is called a \emph{cover}
if $\sum_{i\in C_s} w_i^s \geq  W_s-\lambda_s$, i.e., the set of land parcels in $V_s\setminus C_s$
is not sufficient for satisfying the suitability quota.
Thus, if species $s$ is to be hosted by the reserve, it must hold that
\begin{align}
\sum_{i\in C_s} z_i \geq  u_s,\;\mbox{if $s\in S_1$,\quad\quad or\quad\quad}\sum_{i\in C_s}  x_i \geq u_s,\;\mbox{if $s\in S_2$},\label{eq:cover} \tag{COVER}
\end{align}i.e., at least one land parcel in $C_s$ must be taken.
We note that these inequalities are variants of the well-known \emph{cover}-inequalities for knapsack constraints (see, e.g.,~\citep{KaparisL10}).

The covers $C_s$ can be used to define a stronger version of inequalities~\eqref{eq:zgSC} by replacing $V_s$ with $C_s$ 
(obtaining graphs $G^{C_s}_r$ with arcs $A_s$ replaced by $\{(i,s)\mid i \in C_s\}$). 
The resulting cuts are denoted as \emph{species-cover cuts} (SCC).

\subsection{Constraint-Separation}
\label{subsec:sep}

\paragraph{Separation of Connectivity Cuts and Species(-cover) Cuts}
Let $\rho=(\tilde{\mathbf{u}},\tilde{\mathbf{x}},\tilde{\mathbf{z}},\tilde{\mathbf{y}})$ be the solution of the LP relaxation at the current node of the branch-and-bound tree. 
Connectivity cuts~\eqref{eq:concore},~\eqref{eq:conbuffer} 
(if a model using them is desired), 
species cuts~\eqref{eq:zgSC} and species-cover cuts (SCC) can be separated using maximum-flow computations on a bi-directed graph $D_r$ based on $G_r$ 
(resp., $G^s_r, G^{C_s}_r$), 
where capacities are defined based on $\rho$ 
(see also~\cite{AlvarezLjubicMutzel2013a,AlvarezLjubicMutzel2013b,CarvajalEtAl2011,FischettiEtAl2014}).
While our description applies to~\eqref{eq:concore}, 
the other cuts can be separated in a similar way by a straightforward adaption of the graph used for separation. 
The idea behind the separation is to define a directed graph by splitting the nodes $i \in V$ to directed arcs $(i_1,i_2)$ \red{with associated} capacities $\bar z_i$. 
For $\ell \in V$,
a \red{minimum} $r-\ell$ cut with capacity smaller than $\bar z_\ell$ can then be used to construct a violated inequality~\eqref{eq:concore}.

To be more precise, \red{starting from $G_r=(V \cup \{r\},E \cup A_r)$, the digraph $D_r=(N,A)$ used in the separation of~\eqref{eq:concore} is defined as follows:}
For $i \in V$, define nodes $i_1$ and $i_2$, i.e., $N=\{r\} \cup \{i_1,i_2\mid i \in V\}$ and set 
$A_z= \{ (i_1,i_2)\mid i \in V\}$ and $A'_r= \{ (r,i_1)\mid (r,i) \in A_r\}$.
Bi-direct the edges $\{i,j\}$ in $E$ to arcs $(i,j)$, $(j,i)$. After bi-directing, all ingoing arcs into $i$ are connected to $i_1$, and all outgoing arcs from node $i$ are connected to $i_2$. The set $A$ is defined as the set of arcs obtained this way, plus $A_z$ and $A'_r$, i.e., $A=\{ (i_2,j_1),(i_1,j_2) \mid \{i,j\} \in E\} \cup\{(i_1,i_2)\mid i\in V\} \cup\{(r,i_1)\mid i\in V\}$.
We then define arc-capacities $cap$ as follows based on $\rho$:
\begin{align*}\mathit{cap}_{tv} = \begin{cases}
\tilde z_i, & \text{ if }t=i_1, v=i_2, i \in V, \\
\tilde y_{i}, & \text{ if } t=r, v = i_1, i \in V,\\
\infty, & \text{ otherwise}.
\end{cases}\end{align*}
Thus, arcs $(i_1,i_2)$ obtained by splitting nodes $i \in V$ 
are assigned the value of the corresponding node variable $z_i$, 
root-arcs $(r,i_1)$ are assigned the value of the associated root-variable $y_i$,
and all other arcs (i.e., the ones obtained by bi-directing $E$) 
are assigned infinite capacity. 
It is easy to see that with these capacities, 
any minimum $r-\ell$ cut (for $\ell \in V$) will not involve arcs $(i_2,j_1)$ for $i,j \in V$, as such a cut has infinite capacity. 
Thus the arcs $\bar A_z \cup \bar A_r$ induced by any minimal $r-\ell$ cut are subsets of $A_z \cup A'_r$. 
Any $\bar A_z \cup \bar A'_r$ can be mapped to a set $(\bar W_V,\bar W_A)$ by defining $\bar W_V=\{i\mid (i_1,i_2) \in A_z\}$ and $\bar W_A=\{(r,i)| (r,i_1) \in A'_r \}$.
Now, if
\begin{align*}
\sum_{i\in \bar W_V} \bar z_i  +  \sum_{j\in \bar W_{A}} \bar y_{j} &< \bar z_\ell,
\end{align*}
then $(\bar W_V,\bar W_A)$ defines a violated connectivity cut~\eqref{eq:concore}. 
Thus, the separation problem for~\eqref{eq:concore} and a given $\ell \in V$ can be solved \red{in polynomial time by using a maximum-flow/minimum-cut algorithm.}

In case $\rho$ is integer, the separation procedure can be much simplified:
We construct all connected components $H$ induced by $\bar z_i=1$ 
(by using, e.g., breadth-first search). 
If for such an $H$, \red{$\bar y_i=0$ holds all $i \in H$,} the component is not connected to the root node. A connectivity cut~\eqref{eq:concore} can be constructed by defining $W_A=H$ and $W_V=\{j|\{i,j\} \in E: i \in H, j \not \in H\}$ (i.e., either a node in $H$ must be directly connected to the root node, or one of the nodes neighboring to $H$ must be in the solution) and taking any $\ell \in H$ for the right-hand-side of the cut (we take the node one with smallest index in $H$). 

In both cases (fractional and integer separation), we use the down-lifting by only 
\red{allowing} $j\leq \ell$ for $z_\ell$ on the left hand-side. 
For separation of (SCC), we consider the set $\bar C$ obtained in the separation routine of the cover inequalities as detailed below.

\paragraph{Separating Cover Inequalities}
For a given species $s\in S_1$, a cover $C_s$ can be found by
solving the following separation problem
\begin{align}
\min\left\{\sum_{j\in V_s} \tilde{z}_j q_j  \bigg| \sum_{j\in V_s} w_j^s q_j \geq W_s-\lambda_s\;\mbox{and}\;\mathbf{q}\in\{0,1\}^{|V_s|}  \right\}.\label{eq:coversep}\tag{COVER-Sep}
\end{align}

Clearly, the feasible solutions of~\eqref{eq:coversep} are the covers $C_s$, and the objective function minimizes the sum of the $\tilde z_i$ for the nodes $i$ in the obtained $C_s$. Thus, if the solution of~\eqref{eq:coversep} is smaller than $\tilde u_s$, a violated inequality \eqref{eq:cover} is obtained (for $s\in S_2$, the $z$ need to be replaced by $x$).
The separation problem~\eqref{eq:coversep} is a knapsack-problem in minimization form. We do not solve the problem exactly, but use the following heuristic: First, we sort the nodes in a non-decreasing way by $\tilde z_j/w_j^s$; then, construct we a cover $C_s$ by iteratively picking the nodes sorted in this way, starting with the smallest ratio, until $\sum_{j\in C_s} w_j^s\geq W_s-\lambda_s$.

\paragraph{Implementation of the Cut-Loop}
At any branch-and-bound node, we do not aim to separate all violated inequalities, but follow the outline described in the following list,
and only move from one separation routine to the next, if no violated inequalities
have been found:
\begin{enumerate}
\item Separate~\eqref{eq:cover} and (SCC) or separate~\eqref{eq:zgSC} (depending on the chosen separation strategy, note that (SCC) dominate~\eqref{eq:zgSC}).
\item Separate~\eqref{eq:concore}
\end{enumerate}This scheme is followed in order to avoid excessive calls of the (time consuming) separation routines and also to avoid overloading the LP with too many inequalities.
To further achieve this goal, the separation of~\eqref{eq:concore} is not done for all nodes $\ell \in V$ with $\tilde z_\ell>0$, but only for nodes with $\tilde z_\ell \geq \tau$, where $\tau$ is a given separation threshold parameter (we tried values of $0.5$ and $0.1$ in our computational experiments, see Section~\ref{sec:compres}). 
Moreover, once a violated inequality~\eqref{eq:concore} is found for some $\ell$, we do not consider the nodes $\{i|(r,i) \in W_A\}$ for separation.

For integer solutions of the LP relaxation, only connectivity cuts~\eqref{eq:concore} are separated 
(observe that for correctness of our approach, only this separation is necessary).
In our computational experiments, we also tried \red{other separation strategies, such as separating fractional points only at the root node,  
or separating integer points $\rho$ only.}
Details are given in Section~\ref{sec:compres}.

\subsection{Construction Heuristic and Primal Heuristic}
\label{subsec:heur}

We have designed a construction heuristic, to generate
a starting solution 
for initializing the branch-and-cut, and a primal heuristic, 
which is incorporated inside of the branch-and-cut framework. 
Moreover, to improve the solution found by the construction heuristic, 
we also implemented a local-branching ILP-heuristic~\cite{FischettiLodi2003}.

The use of such heuristics was crucial, since initial computational experiments showed that the internal heuristics of CPLEX, 
the general purpose ILP-solver we used, 
did only find feasible solutions of very poor quality 
(if any solution was found at all). 
This is likely due to the 
symmetric nature of 
our problem,
the structure of the instances and the fact that, 
due to the cutting plane approach, 
only a partial information about the nature of the problem at hand is given to the ILP-solver.

\paragraph{Construction Heuristic: Phase One}
The construction heuristic first 
creates a feasible solution in a greedy fashion, and then tries to remove unnecessary land parcels in a post-processing phase. 
The greedy heuristic is based on the 
TM heuristic by~\cite{takahashi1980approximate} for the Steiner tree problem.
Recall that in the Steiner tree problem, we are given a graph $G=(V,E)$, 
edge costs $c:E\rightarrow \mathbb R_{>0}$ and terminal set $T\subset V$ and we want to find the 
\red{minimum} 
cost tree containing all $T$ terminals.
In the TM heuristic, one starts with a partial solution $\mathcal S$ consisting of a single node from $T$. 
Let $T'=T \setminus \mathcal S$.
Shortest paths to all $t \in T'$ are calculated.
Let $t^* \in T'$ be the terminal in $T'$ with minimum shortest-path distance to $\mathcal S$. 
The terminal $t^*$ and all nodes and edges on the shortest path from $\mathcal S$ to $t^*$ are added to $\mathcal S$. 
This process is repeated, until all terminals are added in $\mathcal S$.

To use a similar heuristic for the \GRSCCB, some adaptations need to be made to deal with the following differences to the Steiner tree problem: 
i) the solution can have up to $k$ components, 
ii) each component consists of a connected core surrounded by a buffer, iii) there is no set of terminals $T$ to be connected, 
but the species protection constraints~\eqref{eq:sizeC} and~\eqref{eq:sizeV} must be fulfilled instead, i.e., 
$P_1$($P_2$) species from $S_1$($S_2$) must be protected in the solution 
(which in turn depends on fulfilling the suitability quota constraints~\eqref{eq:coverZ}, resp.,~\eqref{eq:coverX}). 

In our heuristic (see Algorithm~\ref{alg:heur}), the partial solution $\mathcal S$ is stored as $(\mathcal S_z, \mathcal S_x)$, where $\mathcal S_z$ contains the core nodes of the partial solution, and $\mathcal S_x$ contains all nodes of the partial solution.
For each species $s$, we also keep track of the habitat score of the nodes in the partial solution; these values are stored in $W_s(\mathcal S)$.
The function $u_s(\mathcal S)$ is \texttt{true}, if and only if the suitability quota constraint~\eqref{eq:coverZ}, resp.,~\eqref{eq:coverX} for $s$ is fulfilled by $\mathcal S$ and it is \texttt{false}, otherwise. 
The function $protectedS_1(\mathcal S)$, resp., $protectedS_2(\mathcal S)$ is \texttt{true}, if and only if $P_1$($P_2$) species from $S_1$($S_2$) are protected in the partial solution $\mathcal S$ and it is \texttt{false}, otherwise. 
A solution $\mathcal S$ is feasible if and only if both $protectedS_1(\mathcal S)$ and $protectedS_2(\mathcal S)$ are \texttt{true}. 
The algorithm also uses a function $computeShortestPaths(V',V'')$ which computes the (node-weighted) shortest-paths between $V'$ and any node $j \in V''$ and returns the distances. Moreover, the method  $nodesOnShortestPath(V',V'')$ returns the nodes on the shortest-path between $V'$ and $V''$. More details are given below.

\begin{algorithm}[h!tb]
\DontPrintSemicolon
\KwData{An instance of the \GRSCCB}
\KwResult{A feasible solution $\mathcal S=(\mathcal S_z,\mathcal S_x)$ for the \GRSCCB}
\tcc{the set $T(\mathcal S)$ of terminals is dynamically updated, see the text for details}
randomly pick $k$ nodes from $V$ to initialize $\mathcal S_z$ \;
\red{$\mathcal S_x \gets  \bigcup_{i \in S_z} \delta^+_d(i)$}\;
\While{$\neg protectedS_1(\mathcal S)\lor \neg protectedS_2(\mathcal S)$}
{
$d \gets computeShortestPaths(\mathcal S_z, T(\mathcal S))$\;
$i^* \gets \arg\min_{j \in T(\mathcal S)} d_j$\;
$nodesSP \gets nodesOnShortestPath(S_z,i^*) \cup \{i^*\}$\;
\For{$i \in nodesSP$}
{
$\mathcal S_z \gets \mathcal S_z \cup\{i\}$ \;
\red{$\mathcal S_x \gets \mathcal S_x \cup \delta^+_d(i)$}\;
}
\tcc{update $T(\mathcal S)$}
}
\caption{Construction Heuristic \label{alg:heur}}
\end{algorithm}

During the course of the heuristic, we select nodes $i$ to add to $\mathcal S_z$ based on shortest-path calculations similar to the TM heuristic, i.e., we build connected cores.
We initialize $\mathcal S_z$ by randomly selecting $k$ nodes, 
which ensures that the final solution has at most $k$ connected components. 
Whenever a node is added to $\mathcal S_z$, 
its $d$-neighborhood is added to $\mathcal S_x$ 
(i.e., to accommodate for the the buffer-constraint). 

In contrast to the standard TM heuristic, 
in our approach the terminal set for the shortest-path calculations 
is not known in advance.
Instead, the set of terminals, denoted by $T(\mathcal S)$, 
is dynamically updated in each iteration,
based on the current solution $\mathcal S$ and the associated values of $u_s(\mathcal S)$, $protecedS_1(\mathcal S)$ and $protecedS_2(\mathcal S)$. 
Given a partial solution $\mathcal S$, a node $i$ belongs to $T(\mathcal S)$ if and only if adding it to $\mathcal S_z$ is \emph{helpful} with respect to $protectedS_1(\mathcal S)$ or $protectedS_2(\mathcal S)$.
A node $i$ is deemed \emph{helpful} with respect to $protectedS_1(\mathcal S)$, 
if $protectedS_1(\mathcal S)$ is \texttt{false}, 
and $i \in V_s$ for at least one $s \in S_1$ with $u_s(\mathcal S)=$\texttt{false}, 
i.e., if adding it as a core node increases the habitat suitability score for a species from $S_1$ not hosted by $\mathcal S$, 
and $\mathcal S$ does not already host $P_1$ species from $S_1$. 
Helpfulness of $i$ with respect to $protectedS_2(\mathcal S)$ is defined similarly; 
in addition, all nodes $j \in \delta_d(i)$
are also included in the helpfulness-check (and not just $i$).

To measure the level of helpfulness of a node, we introduce a node-cost function 
$\Delta_i(\mathcal S)$ that is employed
as cost function for the shortest-path calculations (which are done with respect to the node-cost).
The value of the function $\Delta_i(\mathcal S)$ 
is dynamically updated based on the partial solution given by $\mathcal S$.
This is done to take into account that adding a node $i$ to the core 
(i.e., to $\mathcal S_z$) may also cause some other nodes $j$ to be added to $\mathcal S_x$, due to the buffer-constraints. 
Moreover, by $\Delta_i(\mathcal S)$ we also try to take into account  that
the helpfulness of a node $i$ changes,
depending on the suitability scores $w^s_i$ of the node and the species already protected by $\mathcal S$. Let 
\begin{align*}
\mathcal C_i(\mathcal S)=\sum_{j \in \delta^+_d(i) \setminus S_x} c_j,
\end{align*}
\begin{align*}
\mathcal W^1_s(i,\mathcal S)=\Big(w^s_i+W_s(\mathcal S)-\lambda_s\Big)\Big(1-u_s(\mathcal S)\Big)&\quad  \text{for } s \in S_1
\end{align*}
and 
\begin{align*}
\mathcal W^2_s(i,\mathcal S)=\Big(\sum_{j \in \delta^+_d(i) \setminus S_x} w^s_j+W_s(\mathcal S)-\lambda_s\Big)\Big(1-u_s(\mathcal S)\Big)&\quad  \text{for } s \in S_2,
 \end{align*}
where the value \texttt{true} as output of $u_s(\mathcal S)$ is interpreted as one and \texttt{false} as zero.
The value $\mathcal C_i(\mathcal S)$ measures the cost of adding node $i$ to the solution, while $\mathcal W^1_s(\mathcal S)$ and $\mathcal W^2_s(\mathcal S)$ tries to capture the helpfulness of adding $i$ to the solution with respect to a species $s$.
The node cost $\Delta_i(\mathcal S)$ is finally defined as follows:
\begin{align*}
\Delta_i(\mathcal S)=\dfrac{\mathcal C_i(\mathcal S)+0.001}{\Big(\displaystyle\sum_{s \in S_1} \mathcal W^1_s(i,\mathcal S)\Big) \Big(1-protectedS_1(\mathcal S)\Big)+\Big(\displaystyle\sum_{s \in S_2}\mathcal W^2_s(i,\mathcal S)\Big)\Big(1-protectedS_2(\mathcal S)\Big) +0.0001},
\end{align*}
where the value \texttt{true} as output of $protectedS_1(\mathcal S)$, $protectedS_2(\mathcal S)$ is interpreted as one and \texttt{false} as zero.

The construction heuristic is run for $nstarts=20$ different random starting solutions $\mathcal S_z$. 

\paragraph{Primal Heuristic: Phase One}
As a primal heuristic during the branch-and-cut, we use a slightly modified version of the construction heuristic: In the calculation of the node-cost $\mathcal C_i$, we use $c_i(1-\tilde x_i)$ instead of $c_i$. Moreover, the randomly generated starting solutions are constructed by considering the of nodes with $\tilde y_i\geq 0.001$.
In both cases, we run a post-processing procedure, in which we try to remove unnecessary nodes from $\mathcal S_z$, as described below.

\paragraph{Construction and Primal Heuristic: Phase Two (Post-Processing)}
The post-processing is a greedy local improvement procedure in which 
we iterate through the nodes $i \in S_z$, and check, if after removing $i$, 
the solution remains feasible. Note that, together with $i$ we are also removing all nodes from 
${\delta_d}(i)$  which become redundant after removing $i$. Let us denote this set of nodes with 
$\mathcal S^i_x$ where 
$$\mathcal S^i_x = {\delta_d^+}(i) \setminus \bigcup_{j \in S_z, j \neq i} {\delta_d^+}(j).$$
Let $i^*$ be the node whose removal (together with $\mathcal S^i_x$) results in the largest improvement in the objective function. 
We remove $i^*$ from $\mathcal S_z$ and $\mathcal S^i_x$ from $\mathcal S_x$, 
and repeat the process until no additional node can be removed.

Finally,  we point out that a similar heuristic can also be used for \GRSCC\ by setting the buffer size to zero, 
\red{so we have $\delta^+_d(i) = \{i\}$.}

\subsection{Local Branching-based Heuristic}
The solution found by the construction heuristic is further improved using an ILP-based local-search procedure known as \emph{local branching} by\cite{FischettiLodi2003}.
Given a feasible solution $\mathcal S$, local branching explores its $r$-neighborhood by 
employing an ILP-solver in a black-box fashion. 
In the following, we provide specific details of our implementation that deviate from the standard recipe given in~\cite{FischettiLodi2003}, 
following an improved scheme from~\cite{FischettiEtAl2014}.  

In each local search iteration, we start with the basic ILP-formulation of the problem (given in Section~\ref{sec:MIPFormulation}), 
and extend it through an additional local branching constraint which specifies the $r$-neighborhood with respect to $\mathcal S$. In our case, this ILP-formulation is solved through the branch-and-cut, enhanced by the primal heuristic. Even though the complexity  
of the resulting formulation inherits the complexity of the original problem, its feasible region is significantly smaller due to the choice of the parameter $r$. Furthermore, the resulting ILP does not have to be solved to optimality; instead, one interrupts the solver as soon as a feasible solution is found,
i.e., the \emph{first-improvement} local search strategy is applied. In addition, we impose a time-limit for each local search iteration.  If this time-limit is reached, this means that no improving solution is found in the current neighborhood. In the latter case, the size of the neighborhood is increased by $\Delta_r$. 
Whenever a new best solution is found, the size of the neighborhood is reset to $r$. 
The procedure is then repeated with the improved solution $\mathcal S'$ (or with the larger neighborhood), until one of the stopping criteria is satisfied:
(i) the maximum number of local search iterations is reached, (ii) the maximum neighborhood size is reached, or (iii) the overall time limit for the local branching phase is reached.  

We point out that our local branching takes an additional advantage of the primal heuristic as follows: 
if the heuristic manages to produce a feasible solution $\mathcal S'$ which improves upon the currently best known one, but is infeasible with respect to the 
local branching constraint, the current local search iteration is interrupted, and the procedure is repeated by exploring the neighborhood of $\mathcal S'$.

Let $\mathcal S_z$ be the set of $i$ with $z_i=1$ in a given solution $\mathcal S$ and let $r$ be a given radius. The $r$-neighborhood with respect to $\mathcal S$ is defined as a set of all feasible solutions whose Hamming distance with respect to $\mathcal S_z$ is not bigger than $r$. 
Consequently, the following local branching constraint is utilized in our framework:
\begin{align}
 \sum_{i \in \mathcal S_z} z_i \geq |\mathcal S_z|-r. \label{eq:locbra}\tag{LOCBRA}
\end{align}
The constraint~\eqref{eq:locbra} ensures that at least $|\mathcal S_z|-r$ of the core land parcels of the solution $\mathcal S$ also belong to the new solution $\mathcal S'$. 

Upon the termination of the local branching procedure, the branch-and-cut is started. 
The second important and non-standard feature of our local branching implementation is the utilization of a \emph{cutpool}, which collects all violated inequalities detected during the local branching phase. These inequalities are used to initialize the final call of the branch-and-cut procedure. 
The arc-node separator inequalities~\eqref{eq:concore} found during the local branching phase are globally valid, and hence, their collection and recycling through the cutpool significantly influences the overall computing time. Furthermore, the cuts from the cutpool added at the root node also trigger general-purpose cutting planes implemented within an ILP-solver, resulting in stronger bounds at the root node. 
Furthermore, the inequalities~\eqref{eq:concore} from the \emph{cutpool} are also used to initialize each subsequent local search iteration.

In our implementation, we use $r=5$. As a time-limit for the ILPs (i.e., for each single local search iteration) we set 20 seconds.
If no improved solution is found, $r$ is increased by $\Delta_r=5$, until the maximum neighborhood size of $20$ is reached. 
If an improved solution is found, $r$ is reset to five.
These settings have been determined using preliminary computations.

\section{Computational Results}
\label{sec:compres}

In order to assess the effectiveness and suitability of the proposed approach,  we implemented our branch-and-cut framework 
and tested it on three data sets. 
The first data set contains synthetic instances based on grid-graphs, while the second and third data sets are real-life case studies encompassing instances retrieved from the geographic and ecological survey data.

The implementation of the branch-and-cut is done using CPLEX 12.7 as a generic-purpose ILP-solver. All CPLEX-parameters were left at their default values. The experiments were carried out on an Intel Xeon CPU with 2.5 GHz and 16GB of RAM using a single-thread mode.

\subsection{Benchmark Instances}

\paragraph{Grid-Graph Instances}
Following the procedure proposed in~\cite{DilkinaGomes2010,WangOnal2016}, 
our grid-graph instances were generated as follows: A grid of $n \times n$ nodes (set $V$) was created, and an edge $\{i,j\}$ between $i,j \in V$ exists if and only if $i$ and $j$ are adjacent in this grid. 
For each node $i$, the cost $c_i$ was set to an integer value taken uniformly at random from the range $[1,100]$. 
The habitat suitability $w_i^s$ for node $i$ and species $s$ was set to an integer value taken uniformly at random  from the range  $[20,100]$. 
After generating this value, $w_i^s$ is re-set to zero with probability $20\%$ (for $s \in S_1$), resp.\ $10\%$ (for $s \in S_2$) to account for the fact that normally not all land parcels are suited for all species.

We generated four sets of ten instances according to the following scheme:

\begin{itemize}
\item \emph{Set 1}: $n=20$ (hence $V=400$, $E=760$), $|S_1|=1$, $|S_2|=3$.
\item \emph{Set 2}: $n=20$ (hence $V=400$, $E=760$), $|S_1|=3$, $|S_2|=9$.
\item \emph{Set 3}: $n=30$ (hence $V=900$, $E=1740$), $|S_1|=1$, $|S_2|=3$.
\item \emph{Set 4}: $n=30$ (hence $V=900$, $E=1740$), $|S_1|=3$, $|S_2|=9$.
\end{itemize}

We considered the following three conservation scenarios:
\begin{itemize}
\item \emph{Scenario A}: $P_1=|S_1|$, $P_2=|S_2|$.
\item \emph{Scenario B}: $P_1=|S_1|$, $P_2=\lceil 0.5 \cdot |S_2| \rceil $.
\item \emph{Scenario C}: $P_1=|S_1|$, $P_2=0$.
\end{itemize}

In our experiments, the buffer width $d$ is set to one and
scores $w^s_i$ are set to zero for all nodes $i$ at the boundary of an instance, as nodes at the boundary cannot have buffer nodes at one or more sides (as no such nodes exist). 
We considered $k=1$ and $k=3$, i.e., the solution can consist of one connected component, or at most three connected components, 
and we set $\lambda_s=\lceil 0.05 \rceil \sum_{i \in V_s} w^s_i$ for all $s \in S$. 
Thus, in total, we have 240 grid-graph instances (four sets times ten instances times three scenarios times two different values of $k$).

\paragraph{GAP-Instances}
These real-life instances are based on data of the National Gap Analysis Program (GAP), an initiative of the U.S. Geological Survey (USGS)~\citep{GAP,GAPtwo}. 
We considered three states of different size, namely Oregon, Pennsylvania and Vermont. 
In total, 18 real-life instances are generated -- three states times three scenarios (A,B,C, mentioned above) times two different values of $k$.
A detailed case-study which is conducted on this data set is given in Section~\ref{sec:compresCase}.

\subsection{Computational Setting}

In order to analyze the influence of valid inequalities proposed in this paper, and the role of primal and local branching heuristics, 
the following four settings are compared:

 \def\SE{\texttt{OnlyRoot}}
 \def\SF{\texttt{OnlyInt}}
 \def\SG{\texttt{OnlyIntLB}}

\def\SA{\texttt{Basic}}
\def\SB{\texttt{Basic+}}
\def\SC{\texttt{Basic+CP}}
\def\SD{\texttt{Basic+CPLB}}

\begin{itemize}
\item \SA: This is a basic setting in which only arc-node-separator cuts~\eqref{eq:concore} are separated. 
\item \SB: In this setting, in addition to~\eqref{eq:concore},~\eqref{eq:cover} and (SCC) are separated as well.
\item \SC: The setting is the same as \SB, but the construction heuristic and the primal heuristic are turned on.
\item \SD: Finally, in this extension of the \SC\ setting, the local branching procedure has been invoked between the construction heuristic and the branch-and-cut. 
\end{itemize}
By default, in each of the settings, fractional points are separated only at the root node (with at most 20 cuts of type~\eqref{eq:cover}, (SCC) being separated). Since~\eqref{eq:zgSC} are dominated by (SCC), they are not separated in our framework.
The separation threshold $\tau$ for~\eqref{eq:concore} cuts is set to $0.5$. In the experiments for the grid-graphs, we used a time-limit of 1800 seconds. A time-limit of 180 seconds was given to the local branching phase. 

In the following, we provide results of the computational comparison of the four settings on the grid-graph instances, before we provide a detailed analysis concerning the important structural and performance indicators on the set of real-life instances.  

\subsection{Results on Grid-Graph Instances}
\label{sec:compresSyn}

\begin{figure}[h!tb]
\centering
\begin{subfigure}[b]{0.49\textwidth}
  \includegraphics[width=\textwidth] {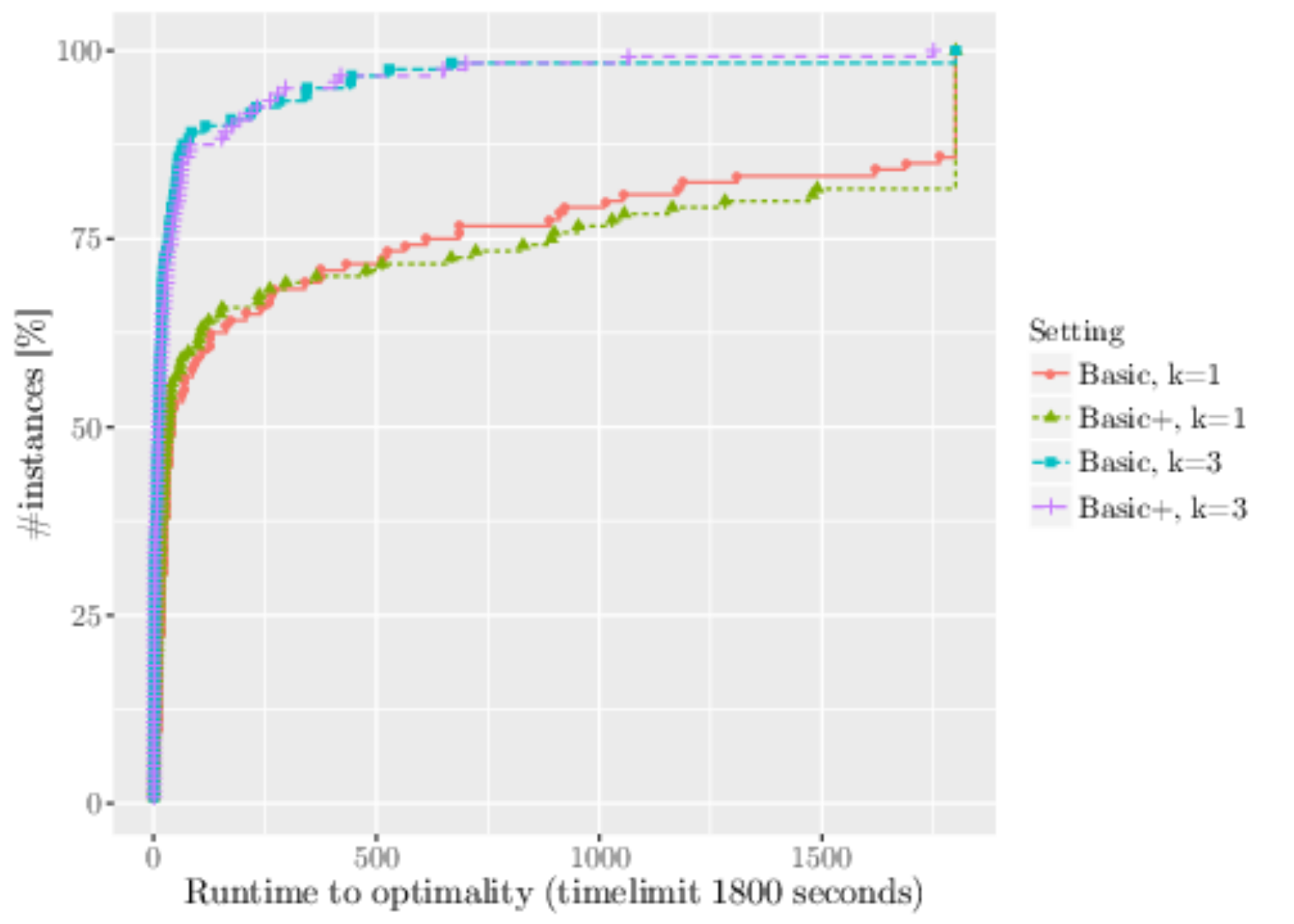}
  \caption{\centering Runtime to optimality\label{fig:runtime}}   
  \end{subfigure}
  \hfill
\begin{subfigure}[b]{0.49\textwidth}
  \includegraphics[width=\textwidth] {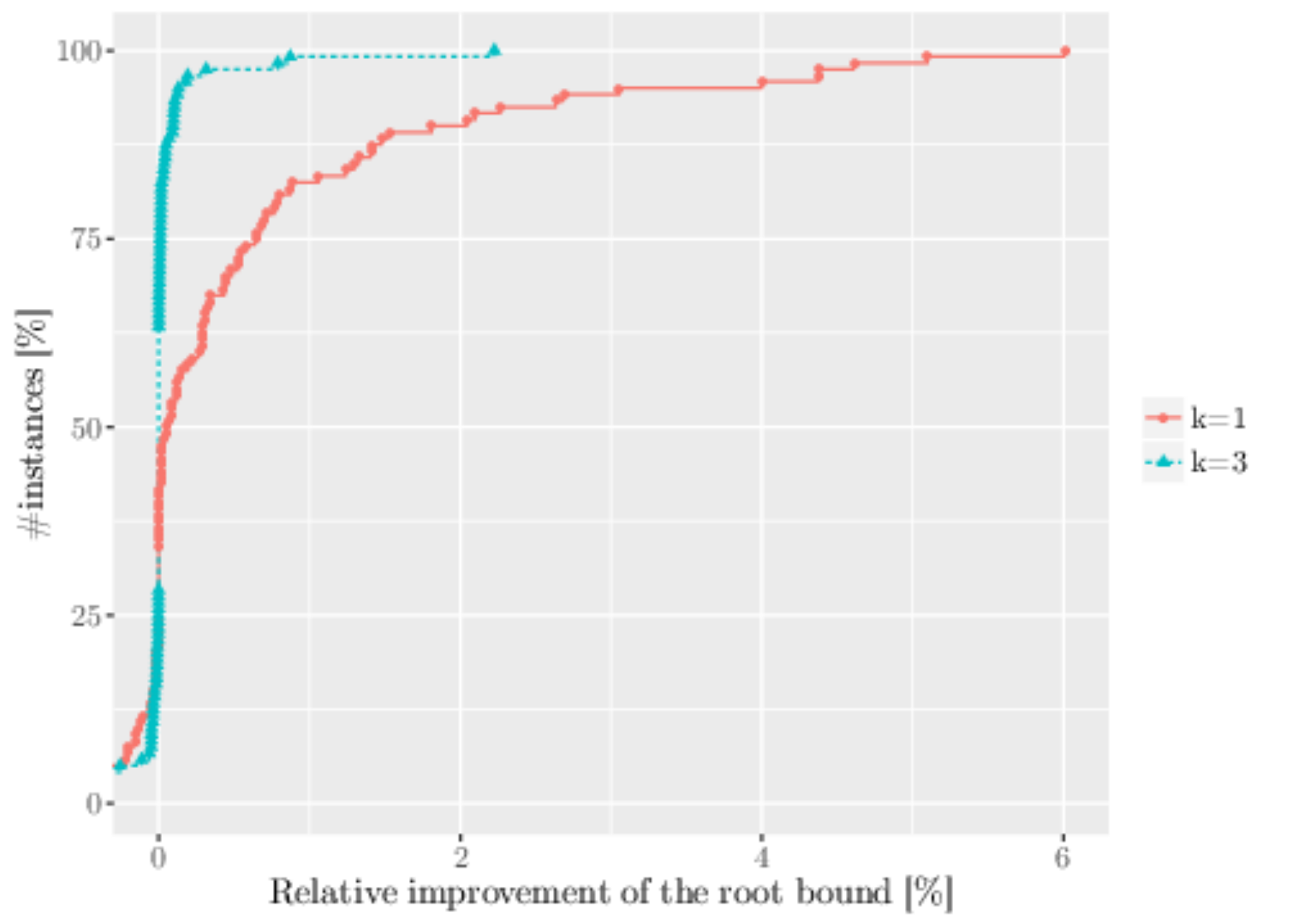}
  \caption{\centering Relative improvement of root node bound\label{fig:rootboundimprovement}}   
 \end{subfigure}
\caption{Comparison of the settings \SA\ and \SB\ for the grid-graphs and $k=1$, $k=3$.}
 \label{fig:SASBcomparison}
\end{figure}

\paragraph{Influence of Valid Inequalities}
In this section we study the influence of the valid inequalities to the quality of bounds obtained at the root node of the branch-and-cut tree (\emph{root bounds} in the following) and the overall computing time. 
We also analyze the problem difficulty with respect to the imposed value of $k$, which is the maximum allowed number of core components. To this end, we compare the computing times to optimality for the settings \SA\ and \SB, and the relative improvement of the root bound obtained by additionally including inequalities~\eqref{eq:cover} and (SCC). 
The results of such comparison are reported in charts of Figure~\ref{fig:SASBcomparison}, which are analyzed below.

The performance chart given in Figure~\ref{fig:runtime} depicts the cumulative computing time for the settings \SA\ and \SB, 
and for $k=1$ and $k=3$. 
In this chart, a point with coordinates $(x,y)$ indicates that for $y$\% of the instances of the considered data set, the total computing time was $\le$ $x$ seconds. 
The first observation that can be made from the obtained result is that solving \GRSCCB\ with a single connected component is computationally much more challenging than solving \GRSCCB\ in which the number of components is relaxed to a greater value. Whereas almost all of the 240 grid-graph instances could be solved to optimality for $k=3$, around 20\% of them remain unsolved within the same time limit if a single core component is imposed. 
Furthermore, for a fixed value of $k$, comparison of computing times between \SA\ and \SB\ reveals that the time invested into a (rather time-consuming) separation of~\eqref{eq:cover} and (SCC) inequalities does not significantly influence the overall performance. This is due to our moderate separation strategy which allows for up to 20 cuts of these types to be added at the root node.  

To compare the influence of the valid inequalities to the quality of the root bound, we report the relative improvement of the root bound obtained from the setting \SB\ with respect to the setting \SA.
Figure~\ref{fig:rootboundimprovement} depicts these values in the cumulative fashion, for $k=1$ and $k=3$. Formally, the relative root bound improvement is defined as $\frac{RB(\SB) - RB(\SA)}{RB(\SA)} \cdot 100\%$, where  $RB(.)$ stands for the the lower bound obtained at the root node of the branch-and-cut tree. Notice that the reported root bounds already take into account the general purpose cuts found by CPLEX, which are turned on by default in all our settings. This also explains why the obtained relative improvements can sometimes take negative values. 

The obtained results indicate that the root bound of the setting \SA\ can be improved by up to 6\%, due to our valid inequalities. 
These relative improvements are particularly pronounced for the more challenging setting in which a single core component is imposed.

We therefore conclude that the computational overhead imposed by the separation of valid inequalities is very moderate compared to the benefits achieved through the improvement of the root bounds, and keep the separation of valid inequalities turned on by default in the remainder of this study. 

\begin{figure}[t!]
\centering
\begin{subfigure}[b]{0.49\textwidth}
  \includegraphics[width=\textwidth] {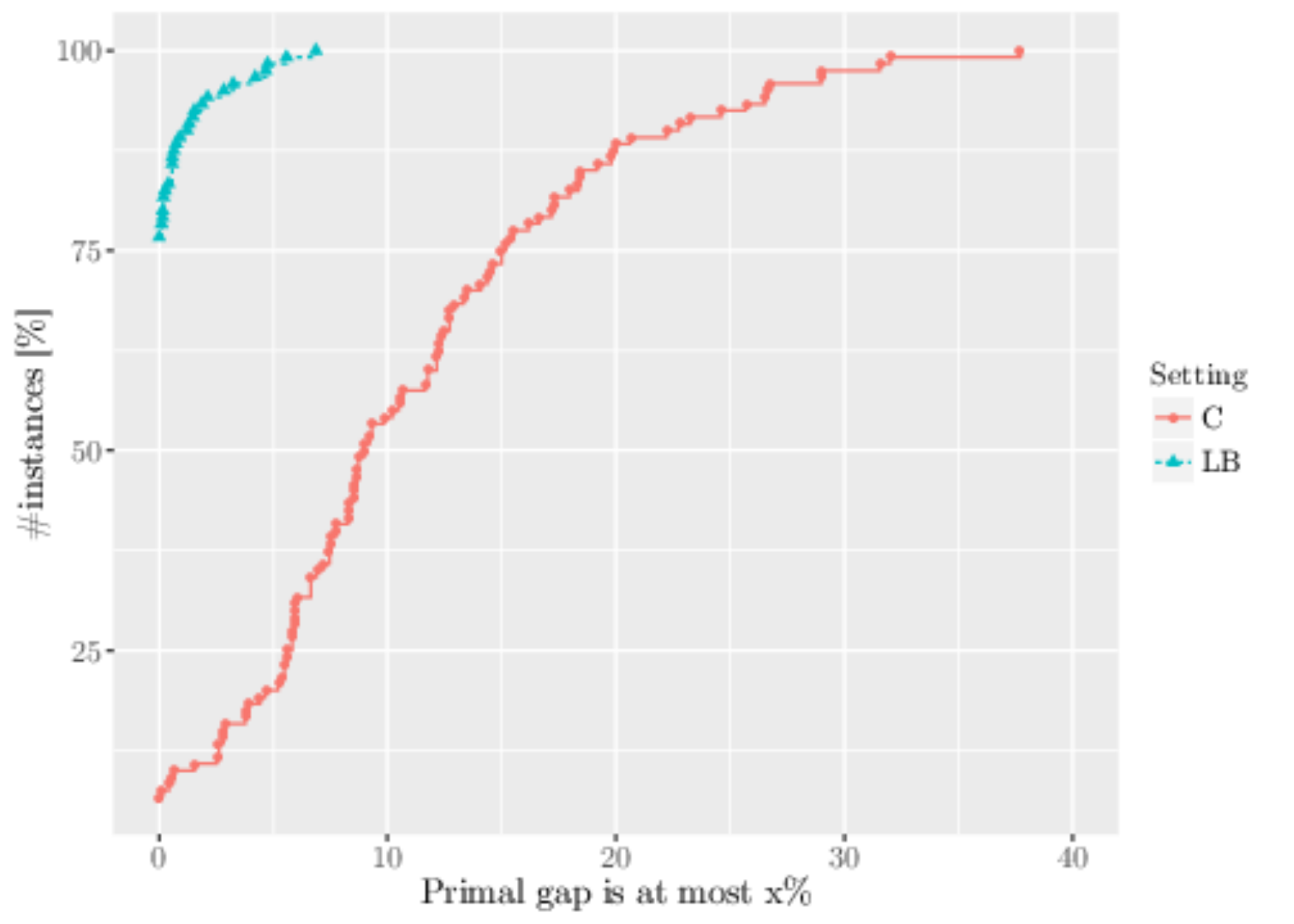}
  \caption{$k=1$}
   \label{fig:heur1}
  \end{subfigure}
  \hfill
\begin{subfigure}[b]{0.49\textwidth}
  \includegraphics[width=\textwidth] {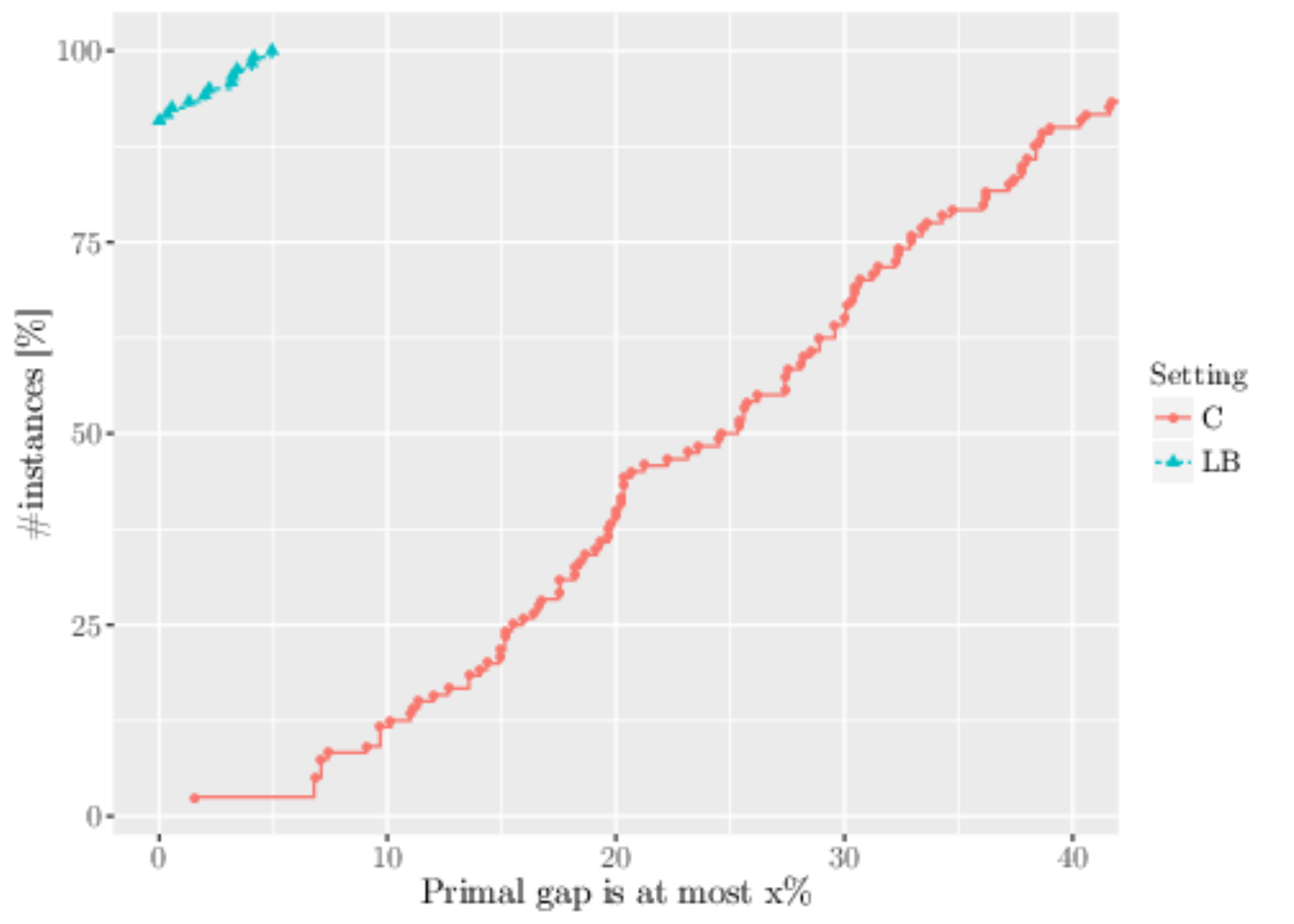}
  \caption{$k=3$}
   \label{fig:heur3}
 \end{subfigure}
\caption{Primal gap of the solutions obtained by the construction heuristic (\texttt{C}) and by local branching (\texttt{LB}) and $k=1$, $k=3$}
 \label{fig:heur}
\end{figure}

\paragraph{Influence of Heuristics}
We now turn our focus on the quality of upper bounds obtained by our solution framework. To this end, we compare the quality of \emph{primal bounds}
used to initialize the branch-and-cut procedure. On the one hand, these bounds are obtained by running the construction heuristic only (setting 
\SC) or, by running the local branching procedure in addition (setting \SD).  
Figure~\ref{fig:heur} analyzes the quality of the construction heuristic and the local branching procedure. In these two cumulative charts, we plot the primal gap $pg[\%]$ against the number of instances: Figure~\ref{fig:heur1} shows the results for $k=1$ and Figure~\ref{fig:heur3} shows the results for $k=3$. 
The primal gap is calculated as $100 \cdot (z^H-z^*)/z^*$, where $z^H$ is the solution value obtained by the construction heuristic/local branching procedure, and $z^*$ is the optimal (or the best known)  solution value for the instance.
The charts show that the construction heuristic works already quite well: for around 90\% ($k=1$), resp., over 60\% ($k=3$) the primal gap is under 20\%. 
Using the local branching procedure after the construction heuristic significantly improves this result, for both $k=1$ and $k=3$.
In case of a single core component, for more than 75\% of the instances, the optimal solution is found upon the termination of the local branching, and the worst primal gap is less than 7\%. Similarly, for $k=3$, local branching finds the optimal solution for about 75\% of the instances, whereas the worst obtained gap remains below 5\%.

We note that there is a moderate computational overhead imposed by the local branching: the average computing time needed for the construction heuristic is below one second, whereas the local branching phase requires around 60 seconds, on average over all 240 grid-graph instances. 
This indicates that local branching does not necessarily pay off for the small instances that can be solved quickly by enumeration in the branch-and-cut tree. On the contrary, high-quality solutions are particularly important for larger instances where the branch-and-cut framework encounters difficulties in closing the final gap. 
We therefore decide to use the setting \SD\ for the case study presented in Section~\ref{sec:compresCase}, since most of the instances considered therein are of the latter type. 

\subsection{Case Study: U.S. Wildlife Conservation}
\label{sec:compresCase}

With the strong increase in the number of endangered species around the globe, several organizations have established sophisticated procedures to gather ecological information of vast areas of the habitat of different types of flora and fauna. Typically,
these habitat surveys are comprised of ecological assessment of the studied area, classification of the species of interest, economical characterizations, and much more. Relevant examples of research groups,
institutes, and public institutions can be found, for instance, in~\citep[][]{ConservationCorr,Microsoft,CBIG}.

A prominent example of an ecology information system is the aforementioned GAP of the USGS. This program follows the methodology proposed in the seminal work of~\citep[][]{ScottEtAl1993}, and it is the result of decades of exhaustive efforts devoted to provide clear, geographically-explicit information on the distribution of native vertebrate species, their habitat preferences, and their management status, in order to determine actual needs in biodiversity protection.

Among the different datasets provided by the GAP program,
it is possible to obtain, for each U.S. state, 
a representation of its territory mapped into so-called hydrological units (HU);
these units and their adjacencies are used to build $G=(V,E)$.
Although the databases supplied by the GAP program have been used before
for analyzing the current wildlife conservation policies (see,~e.g.,~\cite{DrewEtAl2011,LacherWilkerson2014,MeretskyEtAl2012,MinorLookingbill2010}),
we are not aware of OR-oriented papers using this GAP data. 
For this study, we have used data from three U.S. states of varying sizes (Oregon, Pennsylvania and Vermont). For a given state encoded by a graph $G=(V,E)$, 
the problem parameters are obtained as follows. 
\begin{itemize}
\item The set $S$ consists of the terrestrial mammals living in the state according to the GAP data. 
All species which are classified as endangered or vulnerable (at federal or state level) by the corresponding Department of Wildlife
are put into $S_1$ and the remaining species are put into $S_2$. 
The complete lists of animals in set $S_1$ for each considered state is given in Table~\ref{tab:animals}.
As an interesting observation, for Vermont, two of the species (American Marten, Eastern Mountain Lion) listed by the Department of Wildlife do not occur in the GAP-data of the state.

\item In the GAP-data, there are also species distribution models for each species available on a 30 meters $\times$ 30 meters cell basis, i.e., for each such cell and species, there is an "yes/no" flag indicating if the particular cell is suitable for the species. We calculate the score $w_i^s$ of each node $i$ (i.e., each HU) and species $s$ by counting the number of "yes"-cells within the HU.

\item In order to avoid any arbitrary estimation, the cost of taking a land parcel $i\in V$
as part of the reserve is given by the area $a_i$ (in squared meters) of the corresponding HU (see,~e.g.,~\cite{AdamsEtAl2016,JPE:JPE2191,ROMOL-2007}).
\end{itemize}

\begin{table}[t!]
\centering
\caption{Animals in $S_1$ for the considered states} 
\label{tab:animals}
\begingroup\scriptsize
\begin{tabular}{l}
  \toprule
\textbf{Oregon}~(see~\citep{OregonList})  \\
 \midrule
 Canadian Lynx \emph{(Lynx canadensis)}\\
 Gray Wolf \emph{(Canis lupus)} \\
 Columbian White-tailed Deer \emph{(Odocoileus virginianus leucurus)}  \\
 Fisher \emph{(Martes pennanti)} \\
Pygmy Rabbit \emph{(Brachylagus idahoensis)}  \\
 \midrule
 \textbf{Pennsylvania}~(see~\citep{PennList}) \\
 \midrule
 Indiana Bat \emph{(Myotis sodalis)} \\
 Northern long-eared Bat \emph{(Myotis septentrionalis)} \\
 \midrule
 \textbf{Vermont}~ (see~\citep{VermontList}) \\
 \midrule
  Canadian Lynx \emph{(Lynx canadensis)}\\
  Eastern Small-footed Bat \emph{(Myotis leibii)} \\
 Little Brown Bat  \emph{(Myotis lucifugus)} \\
  Northern Bat \emph{(Myotis septentrionalis)}\\
   Indiana Bat \emph{(Myotis sodalis)}\\
   Eastern Pipistrelle \emph{(Pipistrellus subflavus)}\\
   \bottomrule
\end{tabular}
\endgroup
\end{table}

Table~\ref{tab:instances} gives details about the problem instances (column \emph{source} gives the source articles used for classification of species into $S_1$ and $S_2$).

\begin{table}[h!tb]
\centering
\caption{Characteristics of the states under consideration
} 
\label{tab:instances}
\begingroup\scriptsize
\begin{tabular}{l|cccccc}
  \toprule
 $state$ & area (km$^2$) & \#parcels & avg. parcel-area (km$^2$) & $|S_1|$ & $|S_2|$ & source \\
 \midrule
\textbf{Oregon} & 254,799 & 3134 & 81.3 & 5 & 69 & \cite{OregonList} \\
\textbf{Pennsylvania}& 119,280 & 1452 & 82.2 &  2 & 27 & \cite{PennList} \\
\textbf{Vermont}& 24,906  & 301 & 82.7 & 6 & 20 &  \cite{VermontList} \\
   \bottomrule
\end{tabular}
\endgroup
\end{table}

\paragraph{Spatial Analysis and Benefits of the \GRSCCB\ Model}
We now compare the four models addressed in this article for reserve set covering, namely  \GRSC, \GRSCB, \GRSCC\ and \GRSCCB.
Recall that the problems \GRSC\ and \GRSCB\ can be modeled as compact ILP formulations and given to a black-box ILP-solver without any additional interventions. For the models \GRSCC\ and \GRSCCB, the branch-and-cut implementation described in this article is used.

Table~\ref{ta:inforeal} gives a comparison of the results for the real-life instances for \GRSC, \GRSCB, \GRSCC\ and \GRSCCB\ for $k=3$.
The table gives the number of components of the solution ($\#co.$), the number of land parcels ($\#lp.$), the best objective value ($z^*$), and the runtime ($t [s.]$).
The time-limit for these runs was set to 3 hours and an entry $TL$ indicates that the instance could not be solved to optimality within this given time-limit. 
In this case, the number in parentheses next to $TL$ gives the optimality gap,
which is calculated as $100 \cdot (z^*-LB)/z^*$, where $LB$ is the obtained lower bound. 
In addition, Figure~\ref{fig:ratio} gives a plot of $z^\GRSCCB/z^\GRSC$, $z^\GRSCCB/z^\GRSCB$ and $z^\GRSCCB/z^\GRSCC$ (where $z^{\mathcal P}$ is the best solution value obtained for problem $\mathcal P \in \{$\GRSC,\GRSCB,\GRSCC,\GRSCCB$\}$). 
This ratio is an indicator for the potential increase in cost incurred by using the more sophisticated reserve design strategy \GRSCCB\ compared to the simpler strategies 
(of course, a solution of a less-constraint problem may also fulfill the constraints explicitly imposed in \GRSCCB). 
The two compact models (\GRSC, \GRSCB) were given directly to CPLEX, while setting 
\SD
was used for \GRSCC\ and \GRSCCB. 

\red{The obtained results reveal that the model \GRSC\ is the easiest to solve (only one instance remains unsolved within the time-limit, and many are solved within a few seconds only), 
but gives very fragmented reserves, consisting of 12 to 310 connected components. 
The structure of the solutions is slightly better when buffer constraints are imposed: }
Except for OR instance with Scenario A, the solutions for \GRSCB\ consist of at most four connected components, and for five of the nine instances, 
they consist of at most 2 components, i.e., they are even feasible for \GRSCCB\ with $k=3$.
The cost of the solutions for \GRSCC\ and \GRSC\ is very similar, i.e., just imposing connectivity of the reserve only marginally increases the cost of the reserve. 
On the other hand, \red{imposing buffer constraints has a much stronger impact on the solution in terms of the overall cost: the cost for the solutions of \GRSCCB\ (and also \GRSCB) is up to 2.5 times higher than the cost of the solutions without buffer requirements. } 

\begin{landscape}

\begin{table}[h!tb]
\centering
\caption{Number of components, land parcels, solution value and runtime for the best solution for the real-life instances.\label{ta:inforeal}} 
\begingroup\scriptsize
\setlength{\tabcolsep}{1.5pt}
\begin{tabular}{ll|rrrr|rrrr|rrr|rrrr|rrr|rrrr}
  \toprule
   $inst.$ & $sc.$ & \multicolumn{4}{c|}{\GRSC}  & \multicolumn{4}{c|}{\GRSCB} & \multicolumn{3}{c|}{\GRSCC, $k=1$}  & \multicolumn{4}{c|}{\GRSCC, $k=3$} & \multicolumn{3}{c|}{\GRSCCB, $k=1$}  & \multicolumn{4}{c}{\GRSCCB, $k=3$} \\ 
  & & $\#c.$ & $\#lp.$ & $z^*$ & $t [s.] (g[\%])$ & $\#c.$  & $\#lp.$  & $z^*$ & $t [s.] (g[\%])$ & $\#lp.$  & $z^*$ & $t [s.] (g[\%])$ & $\#c.$  & $\#lp.$ & $z^*$ & $t [s.] (g[\%])$ & $\#lp.$ & $z^*$ & $t [s.] (g[\%])$& $\#c.$  & $\#lp.$ & $z^*$ & $t [s.] (g[\%])$ \\ 
 \midrule
VT & A & 25 & 28 & 10333 & 4.40 & 2 & 23 & 22301 & 9.32 & 9 & 11221 & \textit{TL} (8.50) & 3 & 13 & 10343 & \textit{TL} (0.04) & 21 & 23383 & 95.38 & 2 & 23 & 22301 & 56.98 \\ 
  VT & B & 12 & 12 & 9192 & 2.18 & 1 & 19 & 21591 & 4.16 & 7 & 9541 & \textit{TL} (1.13) & 3 & 8 & 9377 & \textit{TL} (0.11) & 19 & 21591 & 37.50 & 1 & 19 & 21591 & 43.71 \\ 
  VT & C & 12 & 13 & 9113 & 0.21 & 1 & 19 & 21591 & 3.80 & 8 & 9391 & 1677.48 & 3 & 8 & 9277 & 795.66 & 19 & 21591 & 29.70 & 1 & 19 & 21591 & 32.28 \\ 
  PA & A & 101 & 124 & 54564 & 2.53 & 5 & 80 & 72559 & \textit{TL} (5.12) & 68 & 55110 & \textit{TL} (0.94) & 3 & 72 & 54623 & \textit{TL} (0.05) & 132 & 117085 & \textit{TL} (85.15) & 3 & 83 & 75630 & \textit{TL} (15.43) \\ 
  PA & B & 26 & 52 & 40936 & 614.68 & 2 & 68 & 56181 & 7649.12 & 52 & 45756 & \textit{TL} (9.11) & 3 & 47 & 43691 & 9135.35 & 64 & 60383 & \textit{TL} (12.83) & 2 & 65 & 56296 & \textit{TL} (3.04) \\ 
  PA & C & 13 & 30 & 22351 & 0.11 & 2 & 56 & 50182 & 71.49 & 38 & 27527 & 5689.42 & 3 & 32 & 22929 & 91.49 & 69 & 58765 & \textit{TL} (13.22) & 2 & 56 & 50182 & 305.19 \\ 
  OR & A & 198 & 237 & 120894 & 300.58 & 7 & 169 & 138558 & \textit{TL} (13.47) & 152 & 122139 & \textit{TL} (0.99) & 3 & 171 & 121688 & \textit{TL} (0.62) & 330 & 284733 & \textit{TL} (133.18) & 3 & 181 & 155624 & \textit{TL} (27.47) \\ 
  OR & B & 43 & 90 & 72369 & \textit{TL} (0.17) & 4 & 151 & 112322 & \textit{TL} (17.34) & 132 & 100445 & \textit{TL} (54.34) & 3 & 100 & 85822 & \textit{TL} (30.20) & 202 & 161255 & \textit{TL} (75.18) & 3 & 140 & 112156 & \textit{TL} (19.40) \\ 
  OR & C & 40 & 73 & 55244 & 4.39 & 4 & 134 & 108024 & \textit{TL} (7.61) & 90 & 68310 & \textit{TL} (18.02) & 3 & 68 & 58181 & \textit{TL} (2.26) & 205 & 162160 & \textit{TL} (72.82) & 3 & 134 & 108035 & \textit{TL} (11.60) \\ 
   \bottomrule
\end{tabular}
\endgroup
\end{table}

\begin{figure}[b!]
\centering
 \includegraphics[width=1\textwidth] {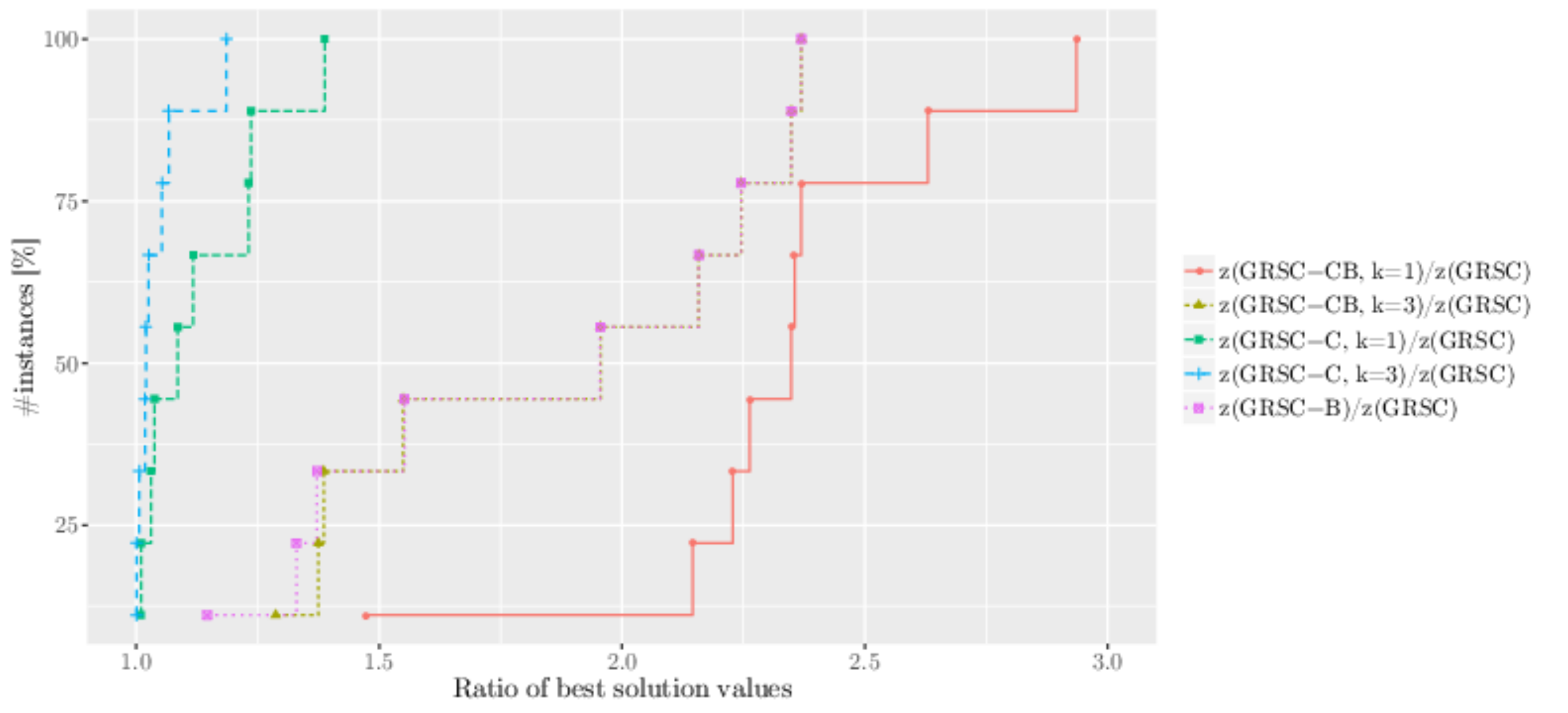}

\caption{Ratio of the best solution values for the different problems}
 \label{fig:ratio}
\end{figure}

\end{landscape}

\begin{figure}[h!tb]
\centering
\begin{subfigure}[b]{0.42\textwidth}
  \includegraphics[width=\textwidth] {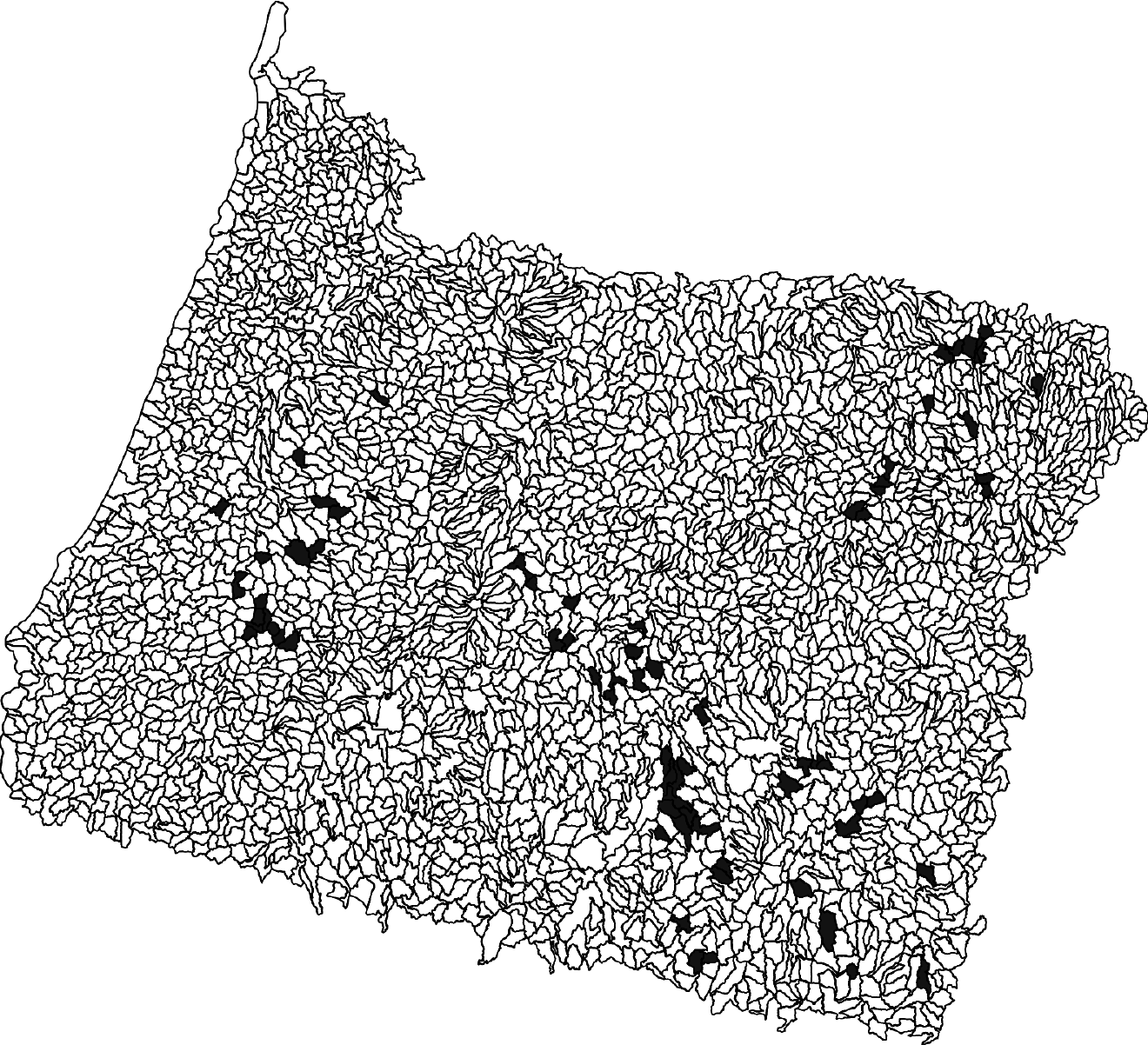}
  \caption{\GRSC \label{fig:OR3GRSC}}
  \end{subfigure}
  \hfill
\begin{subfigure}[b]{0.42\textwidth}
  \includegraphics[width=\textwidth] {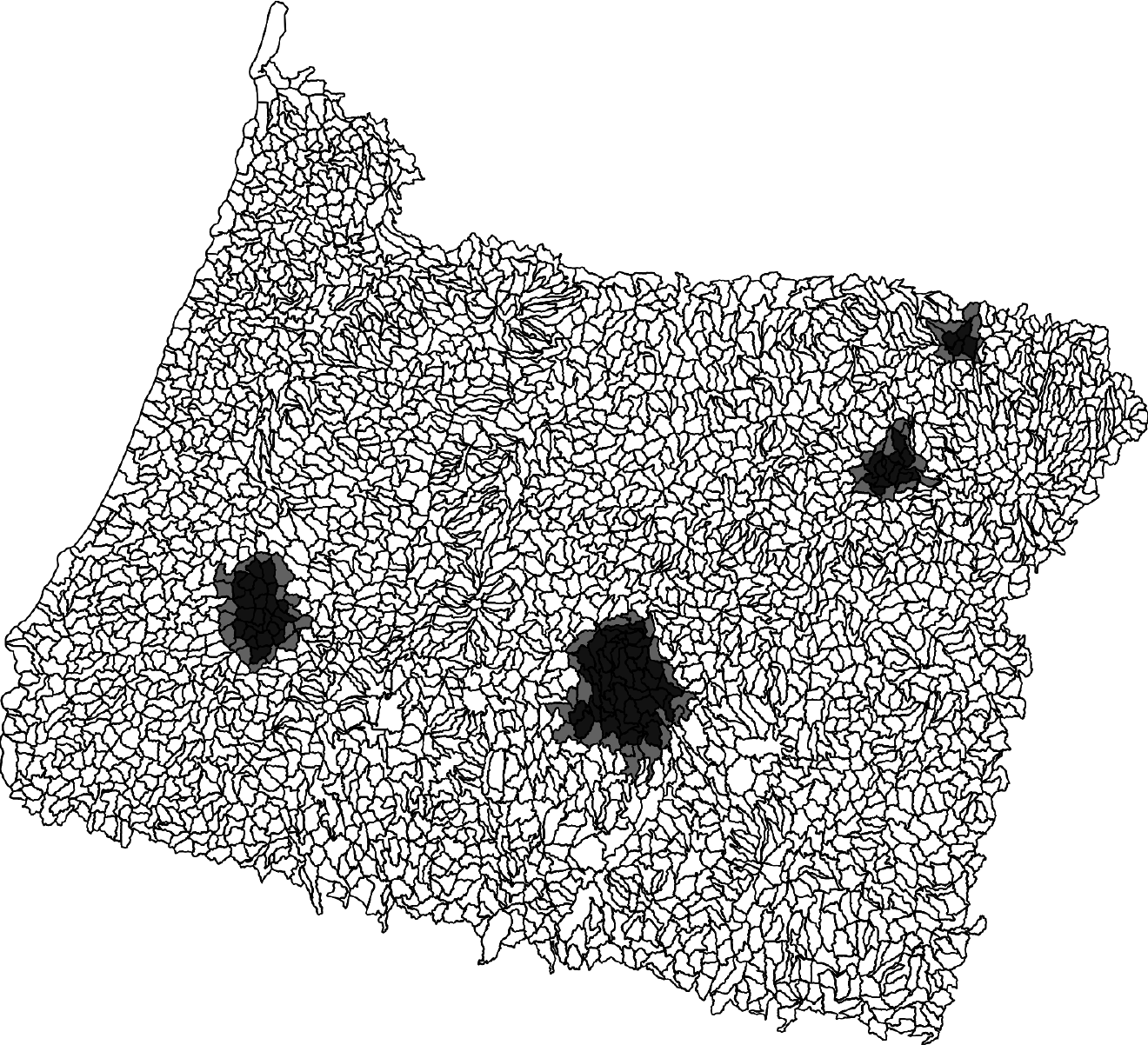}
  \caption{\GRSCB \label{fig:OR3GRSCB}}
 \end{subfigure}
 
 \begin{subfigure}[b]{0.42\textwidth}
   \includegraphics[width=\textwidth] {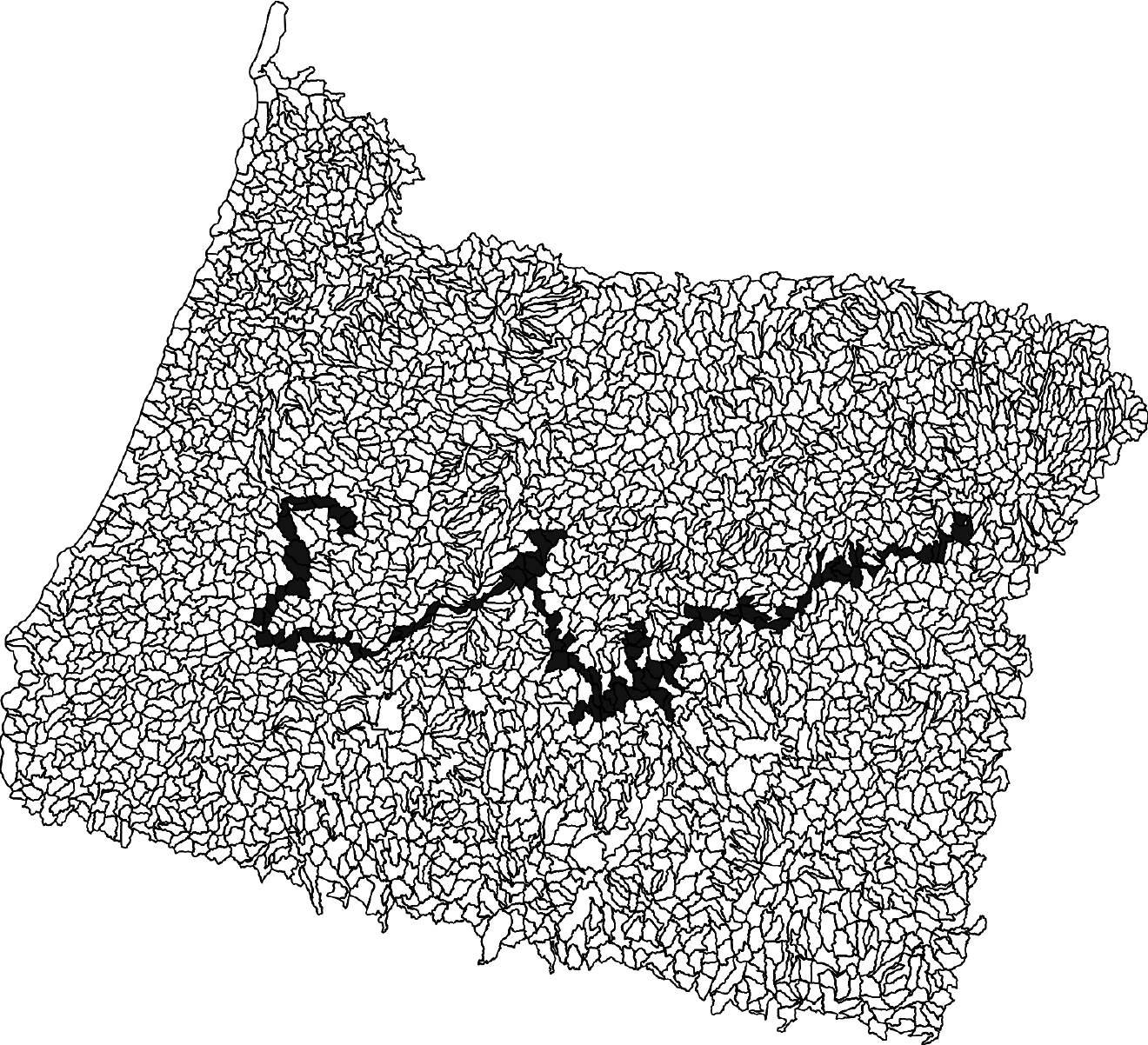}
   \caption{\GRSCC, $k=1$ \label{fig:OR3GRSCC1}}
   \end{subfigure}
   \hfill
 \begin{subfigure}[b]{0.42\textwidth}
   \includegraphics[width=\textwidth] {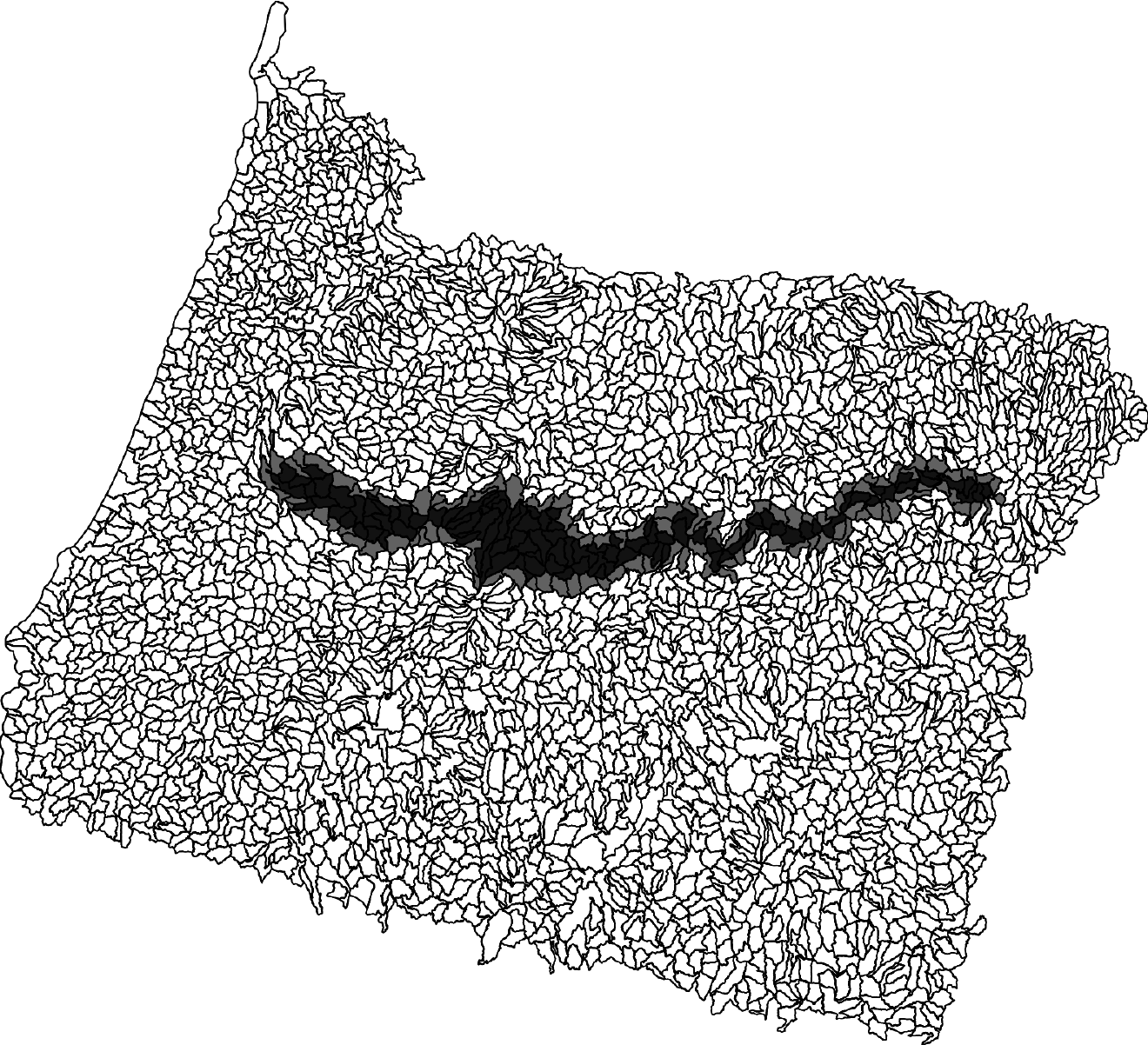}
   \caption{\GRSCCB, $k=1$ \label{fig:OR3GRSCCB1}}
  \end{subfigure}
  
  \begin{subfigure}[b]{0.42\textwidth}
    \includegraphics[width=\textwidth] {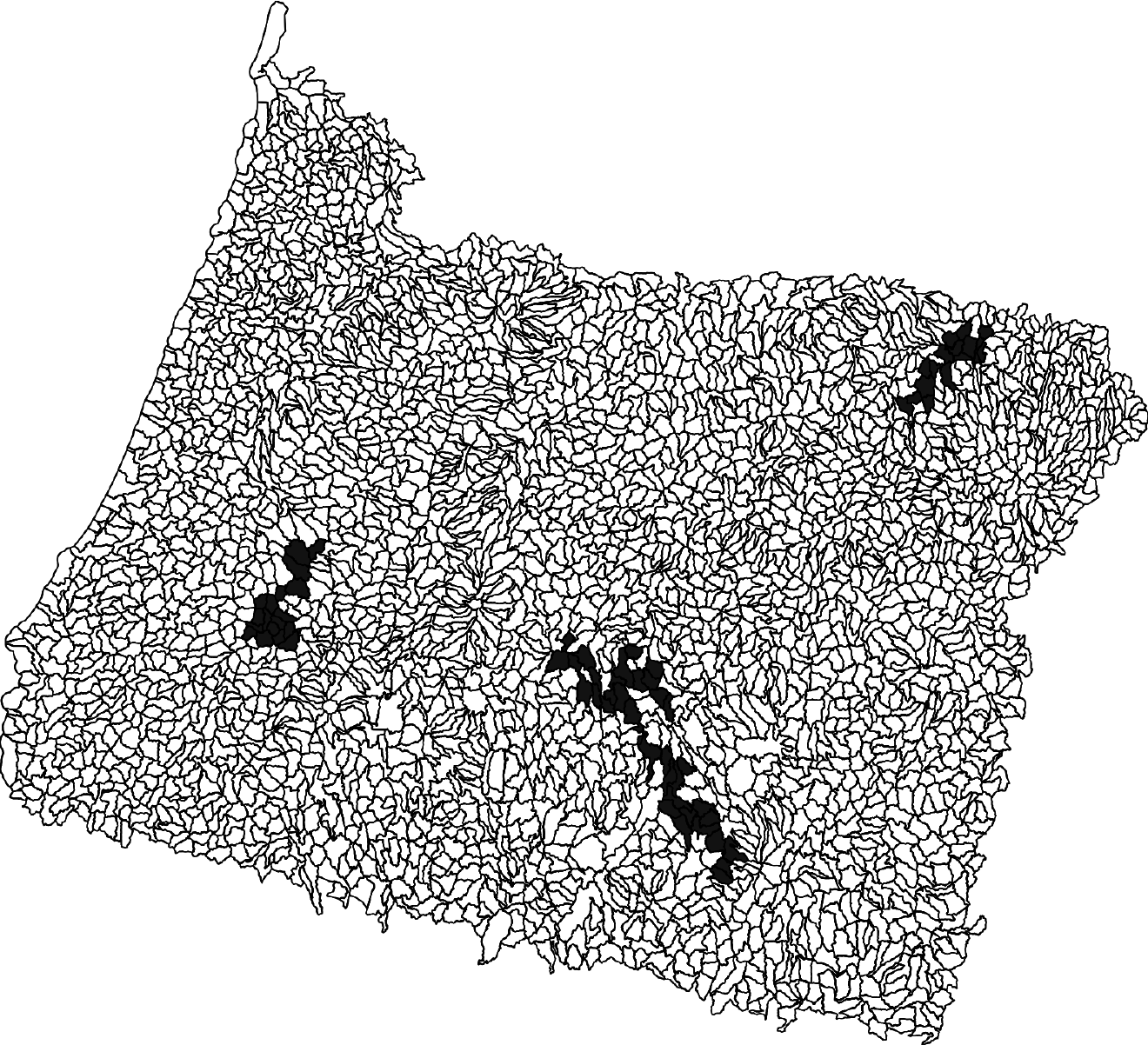}
    \caption{\GRSCC, $k=3$ \label{fig:OR3GRSCC3}}
    \end{subfigure}
    \hfill
  \begin{subfigure}[b]{0.42\textwidth}
    \includegraphics[width=\textwidth] {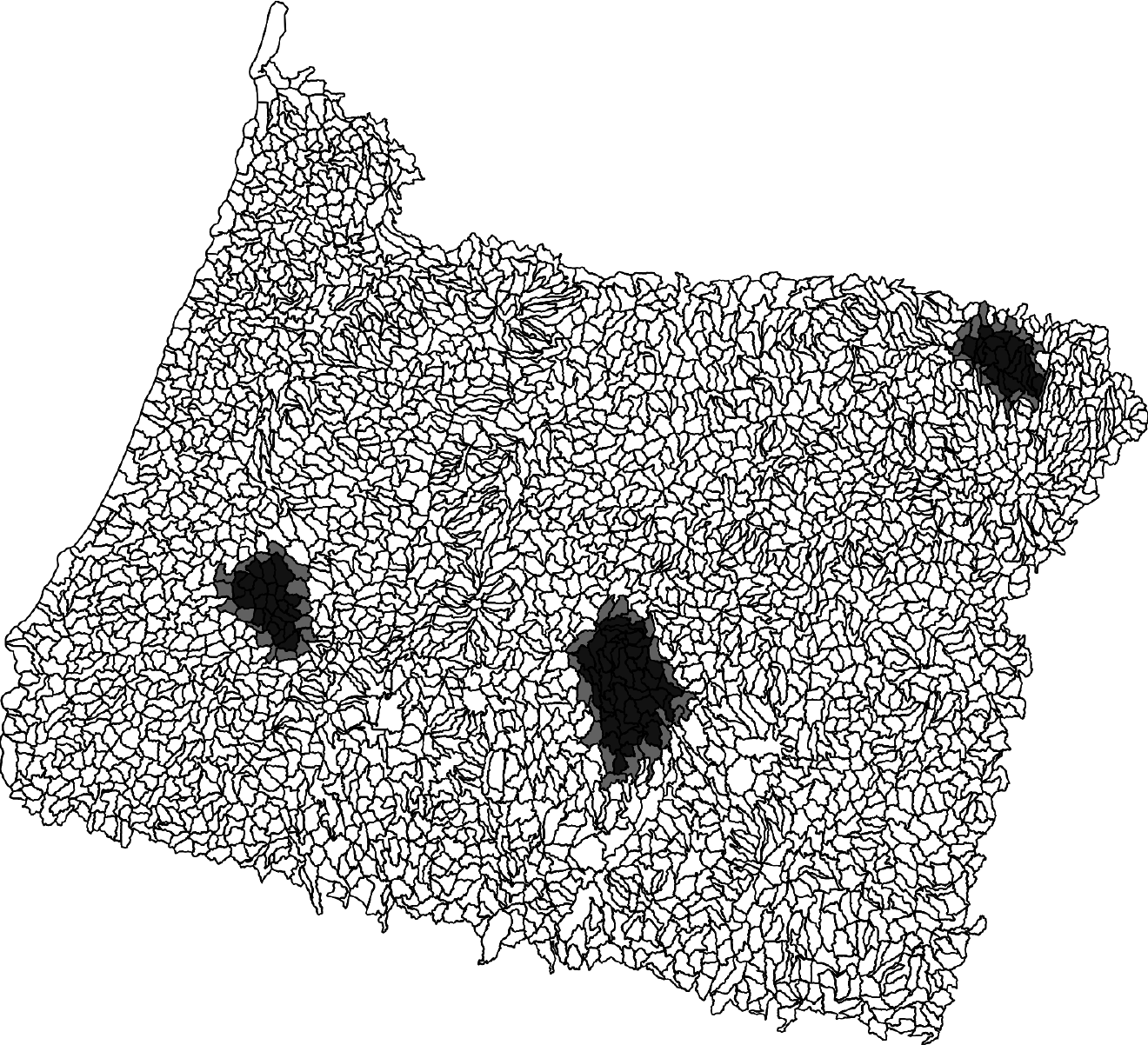}
    \caption{\GRSCCB, $k=3$ \label{fig:OR3GRSCCB3}}
   \end{subfigure}
   
 \caption{Solutions obtained for instance OR and Scenario C \label{fig:OR3}}
\end{figure}

We now look closer into the spacial structure of the obtained solutions, when the four models are applied to the same instance. 
Figure~\ref{fig:OR3} shows the best solutions attained for instance OR and Scenario C for all considered problem variants
(note that for \GRSCC\ and \GRSCCB, we also compare the solutions with $k=1$ and $k=3$). 
The reported solutions show the
potential of the developed framework to design reserves that
respond to different ecological needs. 
It is not surprising that \GRSC\ solutions produce very fragmented reserves as the one depicted in Figure~\ref{fig:OR3GRSC}. 
Such solutions may be suitable when existing protected areas need to be
reinforced by intensifying the protection policy by means of cost-efficient actions.
On the contrary, if the conservation planners seek to design
reserves that are comprised by less scattered units, 
but do not explicitly require the areas to be connected in the strong sense,
the solutions obtained by using \GRSCB\ 
(see e.g., Figure~\ref{fig:OR3GRSCB}) seem quite suitable.
This solution is comprised of only four components, and the components are spatially \emph{compact},
which are both desirable characteristics in many planning contexts.
The solutions reported in Figures~\ref{fig:OR3GRSCC1} and~\ref{fig:OR3GRSCC3} represent optimal solutions for 
\GRSCC\ with $k=1$ and $k=3$, respectively. Comparing the spatial structure of these two solutions, we observe that
different values of $k$ result in spatialy completely different solutions.
This result demonstrates how important for the decision makers is to determine the right choice of the cardinality bound $k$ when reserves with connectivity requirements have to be designed. 
Moreover, the same effect of $k$ holds for the solutions of the \GRSCCB.

Finally, Figures~\ref{fig:OR3GRSCCB1} and~\ref{fig:OR3GRSCCB3} depict the optimal \GRSCCB\ solutions for $k=1$ and $k=3$, respectively. 
Due to the presence of the buffer layer, these solutions are by far more compact than their counterparts  obtained by the alternative three models, \GRSC, \GRSCB\ or \GRSCC. 
This balance of connectivity and compactness makes clear that the \GRSCCB\ formulation is capable of successfully embodying these two ecologically 
functional characteristics to the obtained solutions.

The diversity of the produced solutions are an evidence 
that connectivity and buffer zones are important features that 
greatly influence the spatial layout of the obtained reserves.
Imposing these features allows to define different spatial arrangements,
leading to different ecological profiles that address different
conservation needs. 
For instance, the solutions depicted in Figures~\ref{fig:OR3GRSCC1}
and~\ref{fig:OR3GRSCCB1} resemble wildlife corridors,
which are conservation plans that are suitable
when reserves need be spatially and functionally compatible with other human activities.
\red{Therefore, the \GRSCCB\ model and its variants,
along with the corresponding algorithmic framework given in this paper,
provide a powerful tool for conservation planning. They deliver a flexible decision-aid framework that addresses many different modeling
requirements for designing conservation reserves.}

\clearpage
\subsection{Case Study: Mitchell river catchment}
\label{subsec:mitchell}

As well as U.S., other countries have devoted enormous financial
and human efforts for gathering and analyzing biodiversity information,
and for designing and implementing conservation plans.
One important example corresponds to the Northern Australia Water Futures Assessment program~\cite{NAWFA}, which aims at providing information needed to inform the development and protection of northern Australia's water resources, so that development is ecologically, culturally and economically sustainable.
Within this program, one can find the Northern Australia Aquatic Ecological Assets
initiative~\cite{NAAEA} (NAAE), which focuses on  conducting fine scale assessments of particular catchments in northern Australia.

One of the most important catchment areas approached by the NAAE initiative is the Mitchell river catchment, located in Queensland, northern Australia. 
Within this catchment area, several freshwater fish species were classified as
threatened, and their main threats were spatially and functionally identified.
The studied area (71,630 km$^2$) was divided into 2,316 land parcels (i.e., sub-cachments), and the connections among them corresponded to the river stream network (for further details, see~\cite{CattarinoEtAl2015}).

In the considered region, 46 species (all of them were freshwater fishes) were classified as threatened (the list of species, and their area of presence,
can be found in Table~\ref{subsec:addinfo} in the Appendix).
In Figure~\ref{fig:ausSpecies} we show the spatial distribution density
of the species; darker colors are associated to parcels hosting many species,
while lighter colors are associated to parcels hosting few species.
These species were endangered, due to the presence of four main threats
in the catchment: water buffalo (\textit{Bubalis bubalis}), 
cane toad (\textit{Bufo marinus}), 
river flow alteration (caused by impoundments, channels for water extractions and levee banks),
and grazing land use;
however, a given species was not necessarily menaced by all four threats. 
In Figure~\ref{fig:ausThreats} we show how the number of threats
spatially distributes in the considered area.

\begin{figure}[h!tb]
\begin{center}
\begin{subfigure}[b]{0.49\textwidth}
  \includegraphics[trim={4cm 0 12cm 0},clip,scale =0.21, width=\textwidth] {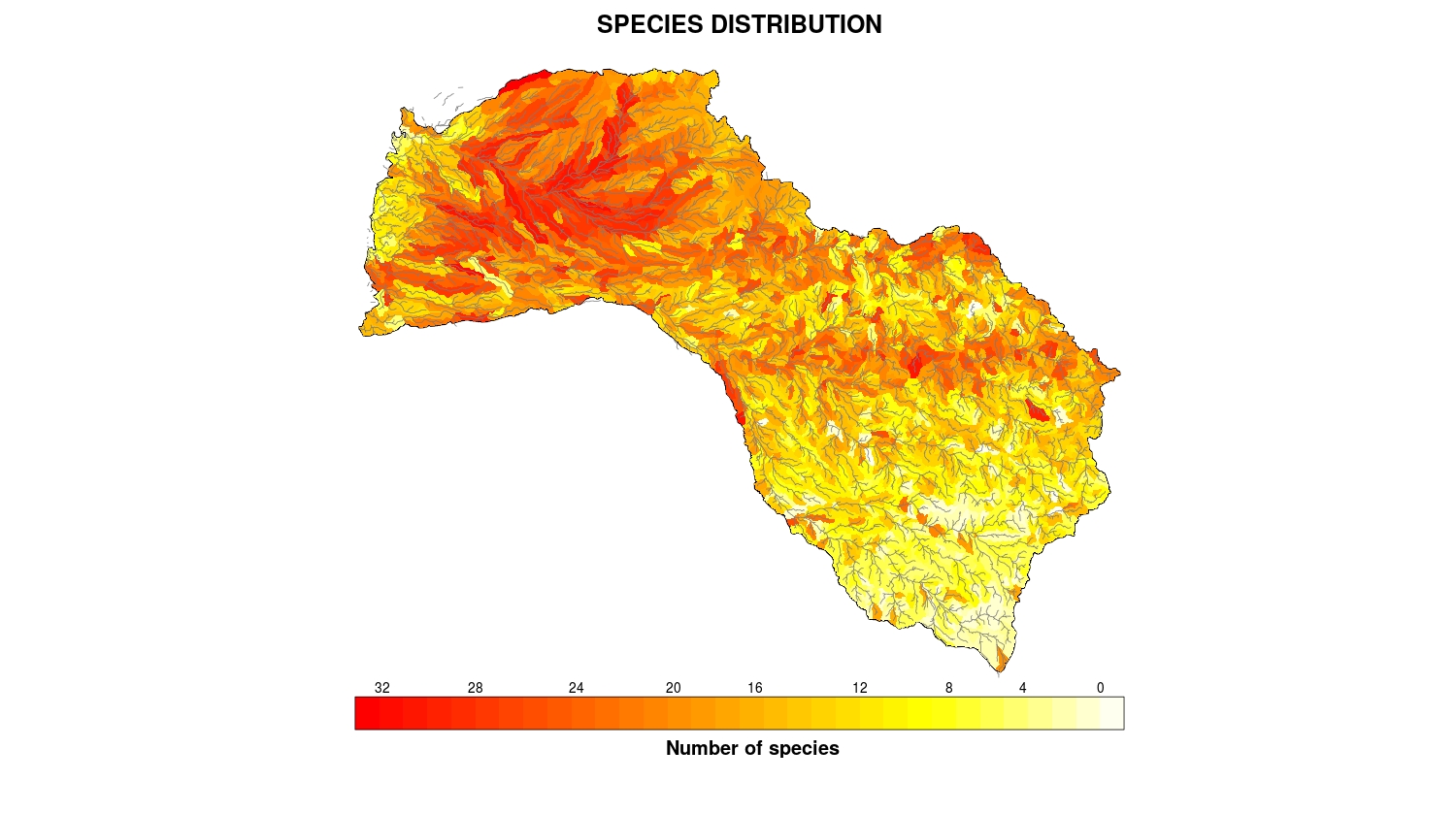}
  \caption{Species distribution \label{fig:ausSpecies}}
  \end{subfigure}
 \hfill
\begin{subfigure}[b]{0.49\textwidth}
  \includegraphics[trim={4cm 0 11cm 0},clip,scale =0.21, width=\textwidth] {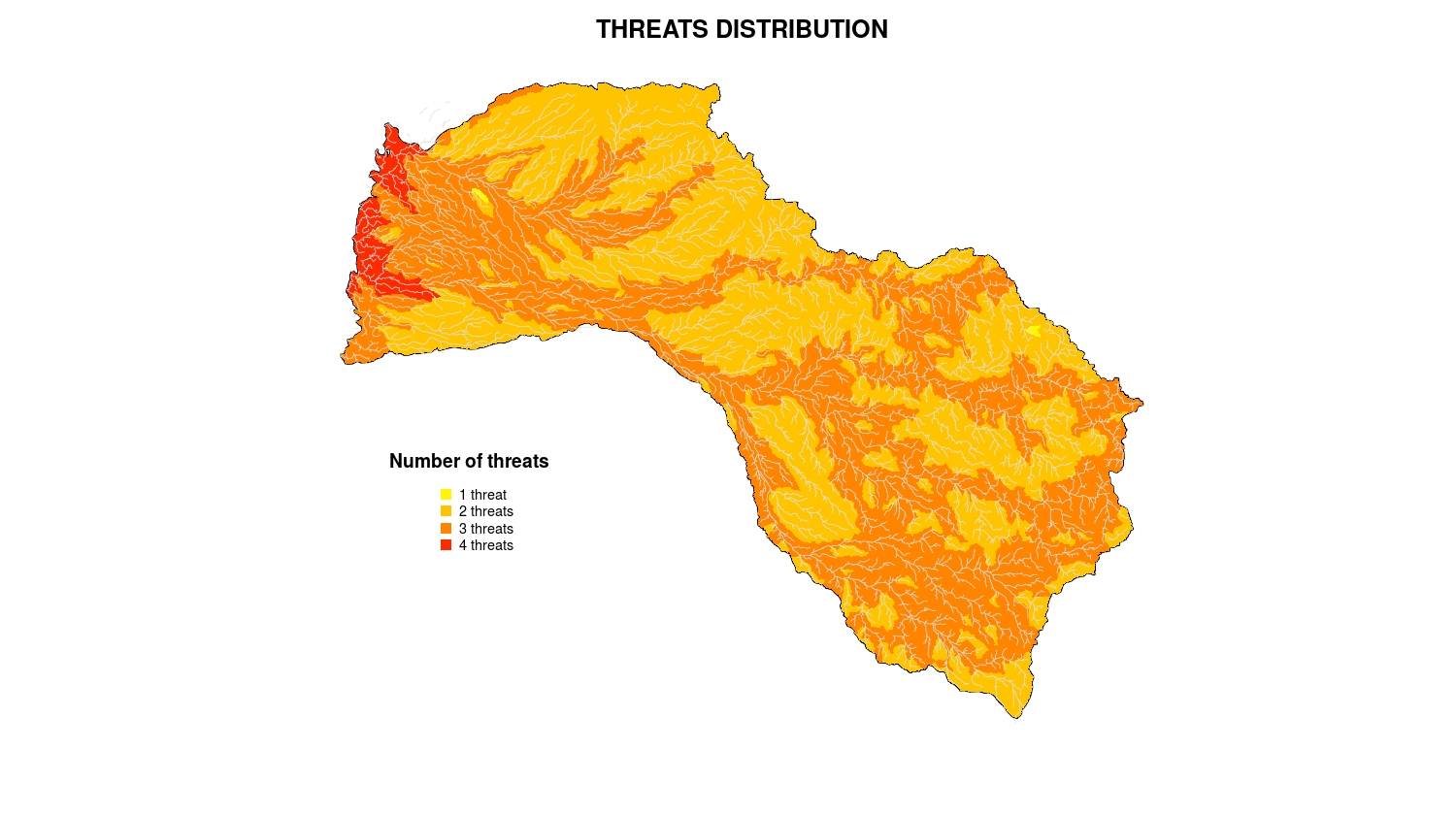}
  \caption{Threats distribution \label{fig:ausThreats}}
 \end{subfigure} 
 \caption{Spatial distribution of species and threats density in the Mitchell river catchment area\label{fig:ausDef}}
 \end{center}
\end{figure}

Since we had explicit information of which species and which threats
co-occurred in each land parcel, the suitability score for a given species
$s\in S$ in a given parcel $i\in I$ is given by
$$
w_i^s = \alpha_i^s \left(1 - \frac{\#threats_i^s}{\#threats_s + 1}  \right)^3 ,
$$where $\#threats_i^s$ corresponds to the number of threats in $i$
that affect $s$, $\#threats_s$ corresponds to total number of threats
that affect species $s$, and $\alpha_i^s$ is a binary parameter taking value 1
if species $s$ is hosted in land parcel $i$, and 0 otherwise.
This expression was defined following the model proposed in~\cite{CattarinoEtAl2015}.
Since in this dataset species were not classified according
to their level of vulnerability, we considered that all of them
were part of $S_1$ (i.e, $|S_1| = 46$).

\paragraph{Obtained results}
For the Mitchell river catchment dataset,
we solved models \GRSC, \GRSCB, \GRSCC\ and \GRSCCB,
using the same algorithmic setting used for the U.S. wildlife dataset.
Since $S_2 = \emptyset$, we imposed the buffer area to be
a \emph{boundary} of thickness given by $d = 1$ (Scenario C).
Due to the structure of the stream network, 
the underlying graph was comprised by several
connected components; this implied that both, the \GRSCC\ and the \GRSCCB\
were infeasible for $k=1$ and $k=2$, 
so we considered $k=3$ in both cases.

\begin{table}[h!]
\centering
\caption{Number of components, land parcels, solution value and runtime for the best solution for the Mitchel river catchment instance.\label{ta:infoaus}} 
\begingroup\scriptsize
\setlength{\tabcolsep}{3pt}
\begin{tabular}{l|rrrr|rrrr|rrrr|rrrr}
  \toprule
   $inst.$ & \multicolumn{4}{c|}{\GRSC}  & \multicolumn{4}{c|}{\GRSCB} & \multicolumn{4}{c|}{\GRSCC, $k=3$} & \multicolumn{4}{c}{\GRSCCB, $k=3$} \\ 
  & $\#c.$ & $\#lp.$ & $z^*$ & $t [s.] (g[\%])$ & $\#c.$  & $\#lp.$  & $z^*$ & $t [s.] (g[\%])$ & $\#c.$  & $\#lp.$ & $z^*$ & $t [s.] (g[\%])$ & $\#c.$  & $\#lp.$ & $z^*$ & $t [s.] (g[\%])$ \\ 
 \midrule
australia & 255 & 402 & 2866 & 0.66 & 5 & 154 & 5527 & \textit{TL} (24.37) & 3 & 192 & 3013 & \textit{TL} (1.28) & 3 & 156 & 5850 & \textit{TL} (36.37) \\ 
   \bottomrule
\end{tabular}
\endgroup
\end{table}

\begin{figure}[h!tb]
\centering
\begin{subfigure}[b]{0.49\textwidth}
  \includegraphics[width=\textwidth] {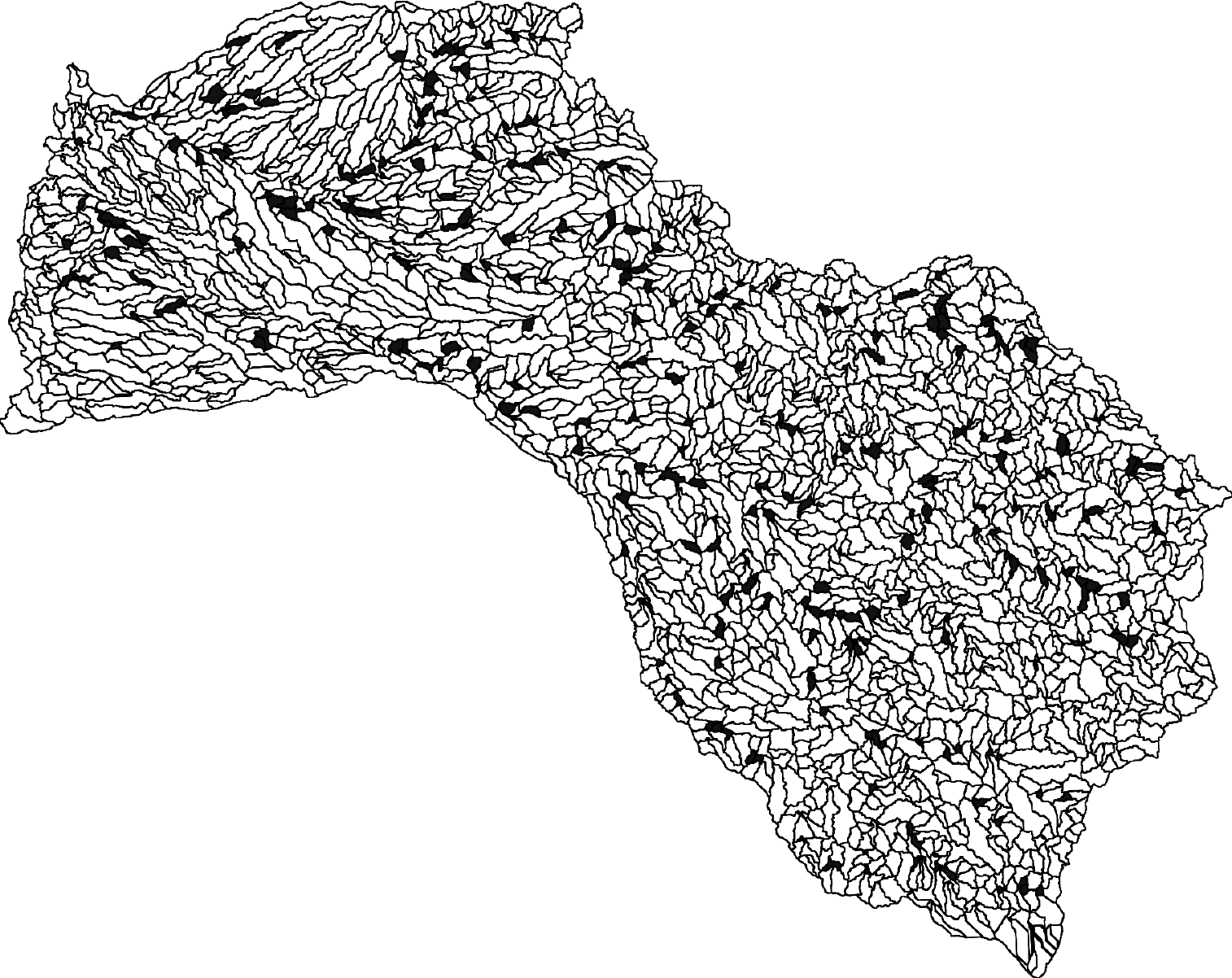}
  \caption{\GRSC \label{fig:ausGRSC}}
  \end{subfigure}
  \hfill
\begin{subfigure}[b]{0.49\textwidth}
  \includegraphics[width=\textwidth] {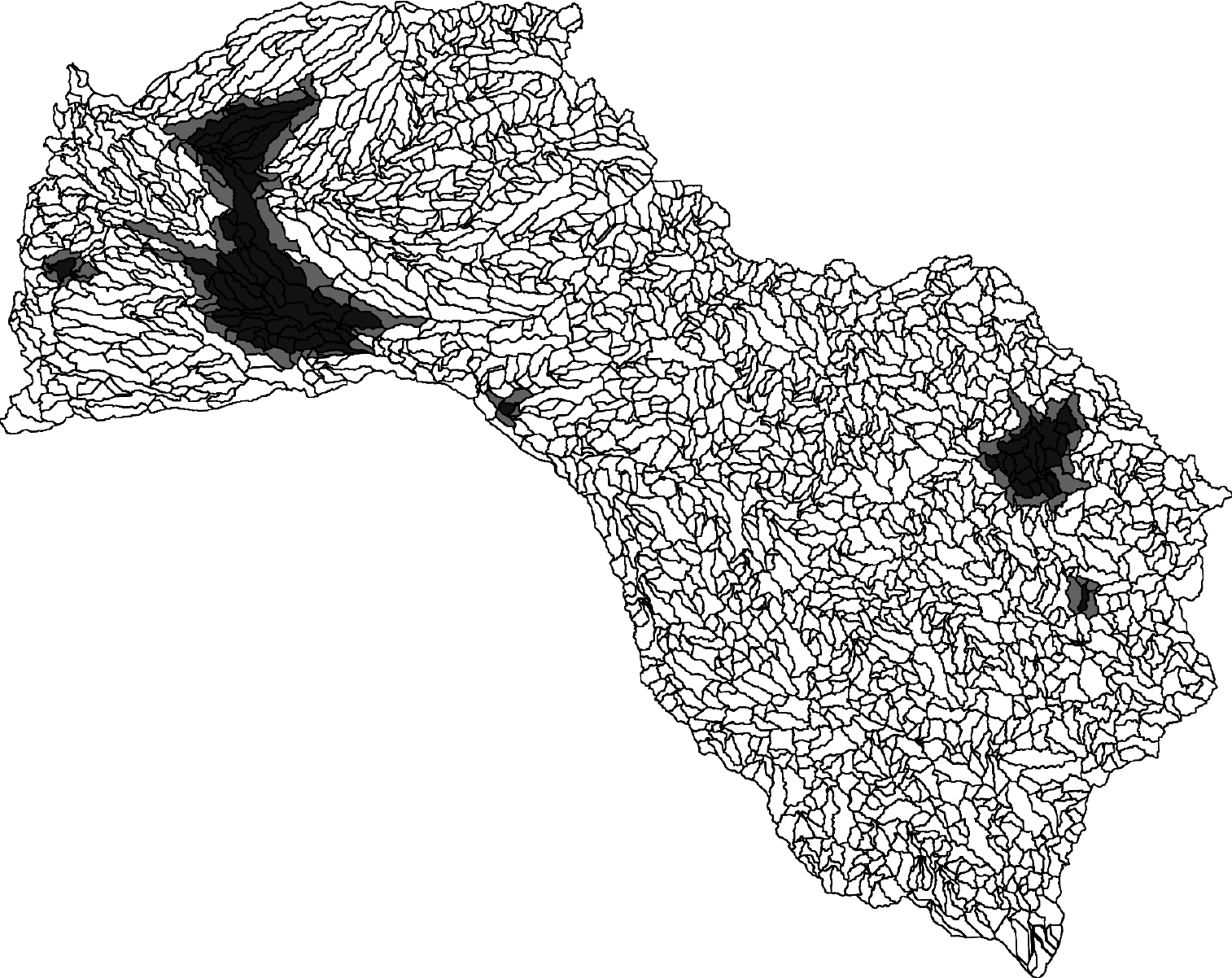}
  \caption{\GRSCB \label{fig:ausGRSCB}}
 \end{subfigure}
  
  \begin{subfigure}[b]{0.49\textwidth}
    \includegraphics[width=\textwidth] {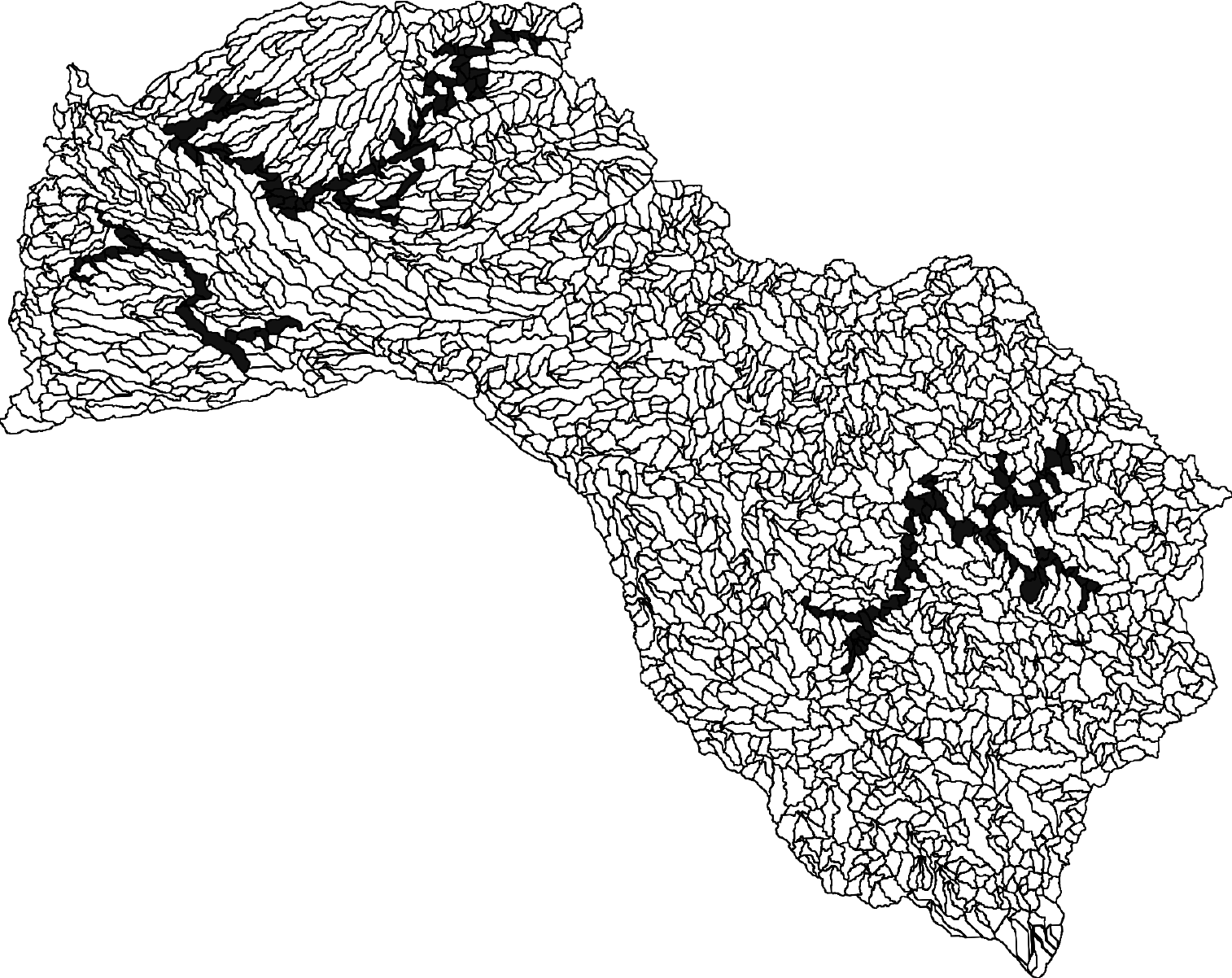}
    \caption{\GRSCC, $k=3$ \label{fig:ausGRSCC3}}
    \end{subfigure}
    \hfill
  \begin{subfigure}[b]{0.49\textwidth}
    \includegraphics[width=\textwidth] {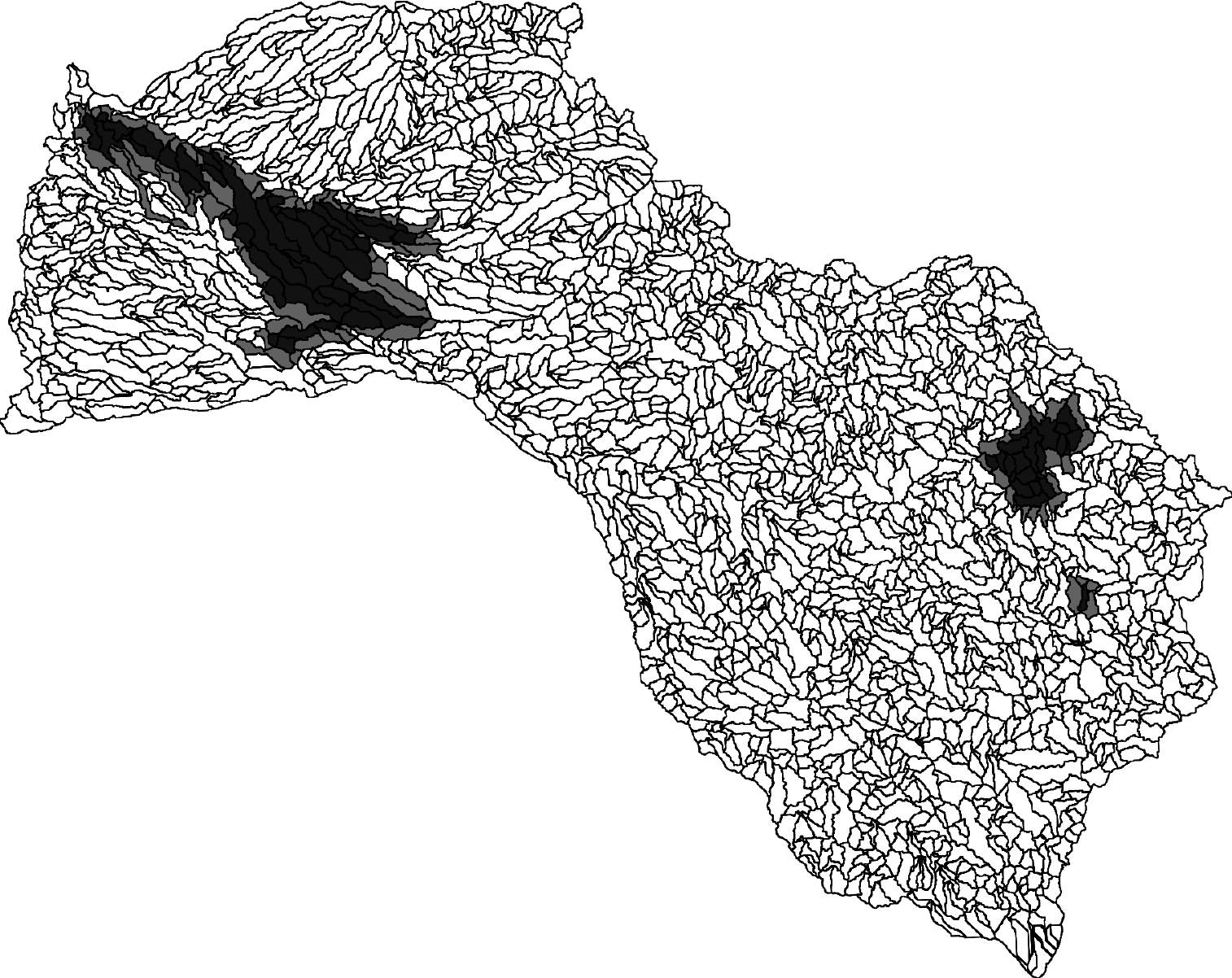}
    \caption{\GRSCCB, $k=3$ \label{fig:ausGRSCCB3}}
   \end{subfigure}
   
 \caption{Solutions obtained for the Mitchell river catchment instance considering Scenario C \label{fig:aus}}
\end{figure}

In Table~\ref{ta:infoaus} (equivalent to Table~\ref{ta:inforeal}) 
we report a summary of the results obtained
when solving the four models on the described instances;
the corresponding solutions are shown in Figure~\ref{fig:aus}.
From the spatial point of view, 
we can observe that the \GRSC\ solution (Figure~\ref{fig:ausGRSC}) is highly fragmented (255 components,
and 402 land parcels),
which is basically due to the presence of many species distributed
along the whole studied area.
When requiring the presence of a buffer (\GRSCB),
the solution changes substantially (Figure~\ref{fig:ausGRSCB}):
the designed reserve is comprised by only 5 components 
(and 154 land parcels in total). 
The solution obtained when solving the \GRSCC\ with $k=3$ (Figure~\ref{fig:ausGRSCC3}),
seems to  properly address the fact that species
live along water flows, encompassing three relatively long
catchment segments. 
Finally, the solution obtained when solving the \GRSCCB\ with $k=3$
(Figure~\ref{fig:ausGRSCCB3}),
spans over similar areas as those spanned by the solution obtained \GRSCC\ with $k=3$,
and it is comprised by almost the same number of land parcels (156 compared to 154). 
From the conservation point of view, the solutions
provided by the \GRSCC\ and the \GRSCCB\ (with $k=3$ in both cases),
allow a more effective implementation of conservation strategies (i.e., implementation of measures against the corresponding threats) 
due to their connectivity and compactness.

From the computational point of view,
it is clear that this dataset is much harder than the ones
considered before.
The two models requiring a buffer boundary
(\GRSCB\ and the \GRSCCB) are considerably more difficult than those that do not require it.
As can be seen from Table~\ref{ta:infoaus},
in these two cases, is not possible to prove optimality within
the running time, and the reached gaps are 24.34\% and 36.67\%, respectively. 
Although the \GRSCC\ optimal solution was not found,
the solution computed within the time-limit reached a 1.28\% gap.
Despite the lack of optimality proof, 
the solutions provided for the different models represent
the first attempt to address from a mathematical programming point
of view, the conservation planning challenges arising in the Mitchell catchment area (see~\cite{Cattarino2016116} for a recent reference on the use of heuristic algorithms for addressing a conservation planning problem in the considered region).

\clearpage

\section{Conclusions}
\label{sec:conclu}

Demographic expansion, natural resource exploitation, and
the consequences of climate change, are among the processes 
that had resulted in a dramatic loss of biodiversity in the last decades. 
The maintenance of biodiversity is crucial for the survival of humankind~\citep[][]{CardinaleEtAl2012}. Thus, immense efforts  have been devoted in the last decades by international organizations, governments, and foundations, for the establishment of protected areas aiming at ensuring
a sustainable landscape for wildlife.
In this paper, we introduced the \emph{Generalized Reserve Set Covering Problem with Connectivity and Buffer Requirements (\GRSCCB)}. This problem is an extension of previous modeling approaches for the design of nature reserves. \red{The problem simultaneously considers, for the first time
\emph{connectivity requirements}, construction of \emph{buffer zones} and \emph{suitability quotes} for species. 
All these constraints } have been identified as being crucial to the design of useful nature reserves. The \GRSCCB\ allows to design a reserve
comprised of one or more connected components; each of them consisting of a core surrounded by a buffer zone, \emph{satisfying minimum suitability requirements for each species.} 
\red{We also consider intermediate problems in which only buffer or connectivity constraints are imposed, denoted by \GRSCB\ and \GRSCC, respectively, and the Generalized Reserve Set Covering Problem (\GRSC) in which both, buffer and connectivity constraints are dropped.
}

We proposed a branch-and-cut framework to solve the \GRSCCB\ and the remaining three problem variants. 
The solution framework is enhanced by the use of valid inequalities and also contains a construction and a primal heuristic, and utilizes a local branching scheme to create feasible solutions of good quality. 
To assess the suitability of our approach, we presented a computational study considering grid graphs, as well as real-life instances representing three different states of the U.S. and a region in northern Australia.
The U.S. instances were constructed using data from the National Gap Analysis Program, an initiative of the U.S. Geological Survey~\citep[][]{GAP,GAPtwo},
with the focus on the protection of mammal species.
The Australia instance  was constructed using
information obtained by the Northern Australia Water Futures Assessment program~\cite{NAWFA}, and corresponds to a conservation setting of
fresh water fishes.

\red{In our study, we compared the solutions obtained by using \GRSCCB\ model against the solutions obtained by 
using less restrictive models (i.e., \GRSC, \GRSCB, and \GRSCC).
On the one hand, we showed the effectiveness of the proposed algorithmic framework on synthetic and real-world instances, by providing optimal or high-quality solutions for many instances of realistic size.
On the other hand, and more importantly, our study
demonstrated the practical versatility of the \GRSCCB\ (and its variants). 
The spatial diversity of solutions produced by using the different problem variants showed that connectivity and buffer zones 
are important features to consider when designing reserves.
Thanks to the flexibility of our models, the developed algorithmic framework provides a powerful tool which allows decision-makers to create desired spatial arrangements while responding to different ecological needs. 
}


\section*{Acknowledgements}
E. \'Alvarez-Miranda acknowledges the support of the National Commission for Scientific and Technological Research CONICYT (FONDECYT grant N.1180670 and Complex Engineering Systems Institute
ICM:P-05-004-F/CONICYT:FB0816), and of the European Union's HORIZON 2020 research and innovation
programme under the Marie Skodowska-Curie Actions (grant agreement N.691149, SuFoRun Project).
The research of M. Sinnl was supported by the Austrian Research Fund (FWF): P 26755-N19 and P 31366-NBL.

\bibliographystyle{plainnat}
\bibliography{Bibliography}{}

\clearpage


\section*{Appendix 1: Additional information of Mitchell river catchment}
\label{subsec:addinfo} \label{sec:appendix}

\begin{table}[h]
{\scriptsize
\centering
\setlength{\tabcolsep}{3pt}
\renewcommand{\arraystretch}{0.75}
\caption{Detailed information of the endangered species included in the Mitchell river catchment area.\label{ta:speciesMitchell}}
\begin{tabular}{lc}
\hline
Name of species					& 	Area (km$^2$) \\
\hline
Scleropages jardinii			&	26130.2		\\	
Nematalosa erebi				&	34153.3\\
Thryssa scratchleyi				&	17161\\
Neoarius berneyi				&	21077\\
Neoarius graeffei				&	8832.2\\
Neoarius leptaspis				&	10920.3\\
Neoarius paucus					&	45154.8\\
Anodontiglanis dahli			&	22921.2\\
Neosilurus ater					&	32947.6\\
Neosilurus hyrtlii				&	26560.3\\
Porochilus rendahli				&	17874.5\\
Arramphus sclerolepis			&	18386.5\\
Zenarchopterus spp.				&	10130.7\\
Strongylura krefftii			&	25112.6\\
Craterocephalus stercusmuscarum	&	54071.5\\
Iriatherina werneri				&	1639.4\\
Melanotaenia splendida inornata	&	70157.5\\
Pseudomugil tennellus			&	2118.9\\
Ophisternon spp.				&	26898.5\\
Ambassis sp.					&	778.9\\
Ambassis agrammus				&	8789.8\\
Ambassis macleayi				&	51412\\
Denariusa bandata				&	11330\\
Lates calcarifer				&	22966.9\\
Amniataba percoides				&	64519\\
Hephaestus carbo				&	10098.4\\
Hephaestus fuliginosus			&	64041.6\\
Variicthys lacustris			&	365.7\\
Leiopotherapon unicolor			&	65926.9\\
Scortum ogilbyi					&	60007.9\\
Glossamia aprion				&	52607.2\\
Toxotes chatareus				&	45386.6\\
Glossogobius aureus				&	40946.1\\
Glossogobius giuris				&	950.8\\
Glossogobius sp. 2				&	24460\\
Hypseleotris compressa			&	370.7\\
Mogurnda mogurnda				&	14594.7\\
Oxyeleotris lineolatus			&	64179.9\\
Oxyeleotris selheimi			&	60793.9\\
Synaptura salinarum				&	3218.8\\
Synaptura selheimi				&	12046.5\\
Megalops cyprinoides			&	10908.8\\
\hline
\end{tabular}
}
\end{table}

\end{document}